\documentclass[journal,transmag]{IEEEtran}
\usepackage{cite, color,soul, amsmath, graphicx, multirow, tabu, cuted} 
\usepackage{flexisym, breqn, cancel}

\ifCLASSOPTIONcompsoc
  \usepackage[caption=false,font=normalsize,labelfont=sf,textfont=sf]{subfig}
\else
  \usepackage[caption=false,font=footnotesize]{subfig}
\fi
\begin{document}
\setlength{\mathindent}{0pt}

\title{The Source Stabilized Galerkin Formulation for Linear \\
Moving Conductor Problems with Edge Elements\textsc{}}

\author{\IEEEauthorblockN{Sujata Bhowmick\IEEEauthorrefmark{$\dag$}
 \textit{Member, IEEE} and  
Sethupathy Subramanian,\IEEEauthorrefmark{$\ddag$} }
\IEEEauthorblockA{ \IEEEauthorrefmark{$\dag$}Department of Electronics System 
Engineering, Indian Institute of Science, Bangalore 560012, India, 
\\ \IEEEauthorrefmark{$\ddag$}Department of Physics and Astronomy, University of Notre Dame, IN 46556, USA}

\thanks{© 2023 IEEE. Personal use of this material is permitted.  Permission from IEEE must be obtained for all other uses, in any current or future media, including reprinting/republishing this material for advertising or promotional purposes, creating new collective works, for resale or redistribution to servers or lists, or reuse of any copyrighted component of this work in other works.}}

% The paper headers
%\markboth{Journal of \LaTeX\ Class Files,~Vol.~14, No.~8, August~2015}%
%{Shell \MakeLowercase{\textit{et al.}}: Bare Demo of IEEEtran.cls for IEEE 
%Transactions on Magnetics Journals}
\markboth{IEEE Magnetics - Accepted article}%
{Shell \MakeLowercase{\textit{et al.}}: Bare Demo of IEEEtran.cls for IEEE 
Transactions on Magnetics Journals}

%-------------------------------------------------------------------------------
%     Part 3.  Abstract, keywords etc.,
%-------------------------------------------------------------------------------
\IEEEtitleabstractindextext{%
\begin{abstract}
The phenomenon of linear motion of conductor in a magnetic field is commonly found in electric machineries such as, electromagnetic brakes, linear induction motor, electromagnetic flowmeter etc. The design and analysis of the same requires an accurate evaluation of induced currents and the associated reaction magnetic fields. The finite element method is a generally employed numerical technique for this purpose. However, it needs stabilization techniques to provide an accurate solution. In this work, such a stabilization technique is developed for the edge elements. 
The stability and hence the accuracy is brought in by a suitable representation of the source term. The stability and accuracy of the proposed scheme is first shown analytically and then demonstrated with the help of 2D and 3D simulations. 
The proposed scheme is parameter-free and it would require a graded regular mesh along the direction 
of motion.

%fluid dynamics

\end{abstract}

% Note that keywords are not normally used for peerreview papers.
}

\maketitle

\IEEEdisplaynontitleabstractindextext

\IEEEpeerreviewmaketitle

%-------------------------------------------------------------------------------
%     Part 4. Introduction 
%-------------------------------------------------------------------------------
\section{Introduction}
\IEEEPARstart{T}{he} numerical simulation of electrical machineries and equipment is inevitable for their economical design and safe operation. The finite element method (FEM) is a commonly employed numerical technique. The FEM is known to produce highly accurate solutions for second order diffusive simulations. The same is not true when dominant first order terms are present. The governing equations of conductor moving in a magnetic field fall into this category. Consider the following governing equations of conductor moving in a magnetic field \cite{reducedb, fmbase},
\begin{equation} \label{eqge1}
\sigma \nabla \phi ~-~  (\nabla \cdot \dfrac{1}{\mu} \nabla) {\bf{A}} - \sigma~ {\bf{u}} \times \nabla \times {\bf{A}} = \sigma~ {\bf{u}} \times {\bf{B_{a}}}
\end{equation}	
\begin{equation} \label{eqge2}
\begin{split}
\nabla \cdot (\sigma \nabla \phi) - \nabla \cdot (\sigma ~ {\bf{u}} \times \nabla \times {\bf{A}}) =  \nabla \cdot (\sigma ~ \bf{u} \times \bf{B_{a}})
\end{split}
\end{equation}	
where, $\phi$ is the scalar potential arising out of the current flow, 
$\bf A$ is the magnetic vector potential associated with reaction magnetic field 
$\bf b$,  $\bf u$ is the velocity of the moving conductor, $\mu$  is the 
magnetic permeability and $\sigma$ is electrical conductivity.

It can be seen that, for the variables $\bf{A}$ and $\phi$, all the derivatives in equation (\ref{eqge2}) are second derivatives. So, this equation is not expected to introduce any instability in the solution. However, the same is not true for the first equation (\ref{eqge1}); here, the first derivative is present in the form of $\sigma~ {\bf{u}} \times \nabla \times {\bf{A}}$. When this becomes dominant, more precisely, when the quantity $\mu\sigma{\bf{|u|}}\Delta z/2$ becomes larger than 1, the numerical instability ensues; { where $\Delta z$ is the element length along the flow direction} \cite{cdbook, upfdm1_sp, up1}.
This quantity is called as Peclet number ($Pe = \mu\sigma{\bf{|u|}}\Delta z/2$).  

In such a situation, to bring in stability and accuracy to numerical solutions, several numerical remedies have been proposed. Among these, the upwinding techniques are commonly used across disciplines. The upwinding schemes are proposed for the fluid dynamics transport equation and extended for the moving conductor problems \cite{mc6tf1, mc6tf2, mc2av1, mc3eb1, mc4ge1, mc5mc1,  mc2av2, mc3eb2, mcsupg_mcfit, mcsupg_cable, mcsupg_mfluid}. 
The upwinding-based schemes can be inferred as to bring in the stability by  introducing the right amount of diffusion \cite{quada1,revfic1}. The \emph{correct} amount of diffusion is decided by the stabilization parameter $\tau$, which is defined to be, $\tau = \coth(Pe) - 1/Pe$.
On the other hand, the recent source-stabilized finite element schemes are primarily proposed for the linear moving conductor problems \cite{su1, su2}. They do not seek stability  by adding diffusion (upwinding) to the governing equation. Instead, stability is
brought in by the appropriate representation of the source term, which
mitigates numerical instability via pole-zero cancelation. In addition, the source-stabilizing schemes are shown to be free of non-physical currents at the material boundary \cite{sus3}. 

It can be noted that, all of the above mentioned stabilisation techniques are derived for the linear nodal elements for the one-dimensional problem of equal discretisation, and they are heuristically extended for the 2D and 3D problems \cite{quada1}. The one exception can be the source-stabilized scheme proposed in \cite{su2}, where stability is analytically shown for a simplified 2D problem. 

In electrical engineering, edge elements are widely used to accommodate for the discontinuity of the normal field at the material interfaces; this is not possible with the nodal elements. 
It can be noted that, numerical instability at high velocities is present in the edge element formulation as well. In order to cater this, there are upwinding techniques proposed for the moving conductor problems with the edge elements \cite{edgeup2, r:edgeup3}. These are generally based on the heuristic extension of the upwinding techniques proposed for the fluid dynamics transport equation. Hence, they are also susceptible to transverse-boundary error at the material interfaces \cite{ss1, soldreview1, soldreview2, nemosu}.

In this work, an attempt is made to propose a source-stabilized Galerkin finite element formulation for the edge elements. For this, a simplified version of the moving conductor problem is considered; using that, the stability of the proposed formulation is established. Then, in order to correctly represent the edge elements, as well as, the curl nature of the governing equation, an extensive stability analysis is carried out in 2D. Subsequently, numerical exercises are carried out both in 2D as well as 3D. 

In the next section, description of the present work is provided and it starts with the stability analysis for a simplified problem.

%Then, considering the coupled curl nature of the governing equations, the numerical stability of the proposed method is shown for a simplified 2D problem. Subsequently, more practical 2D and 3D problems are taken up for the numerical simulations and the stabilisation of the proposed method is demonstrated.     

%-------------------------------------------------------------------------------
%     Part 5. Present Work
%-------------------------------------------------------------------------------
\section{Present Work}
%-------------------------------------------------------------------------------
%     Part 5.1. 1D Analysis
%-------------------------------------------------------------------------------
\subsection{Analysis with limiting 1D version of the problem} \label{sec1D.A}

\begin{figure}
		\centering
		\includegraphics[scale=0.33]{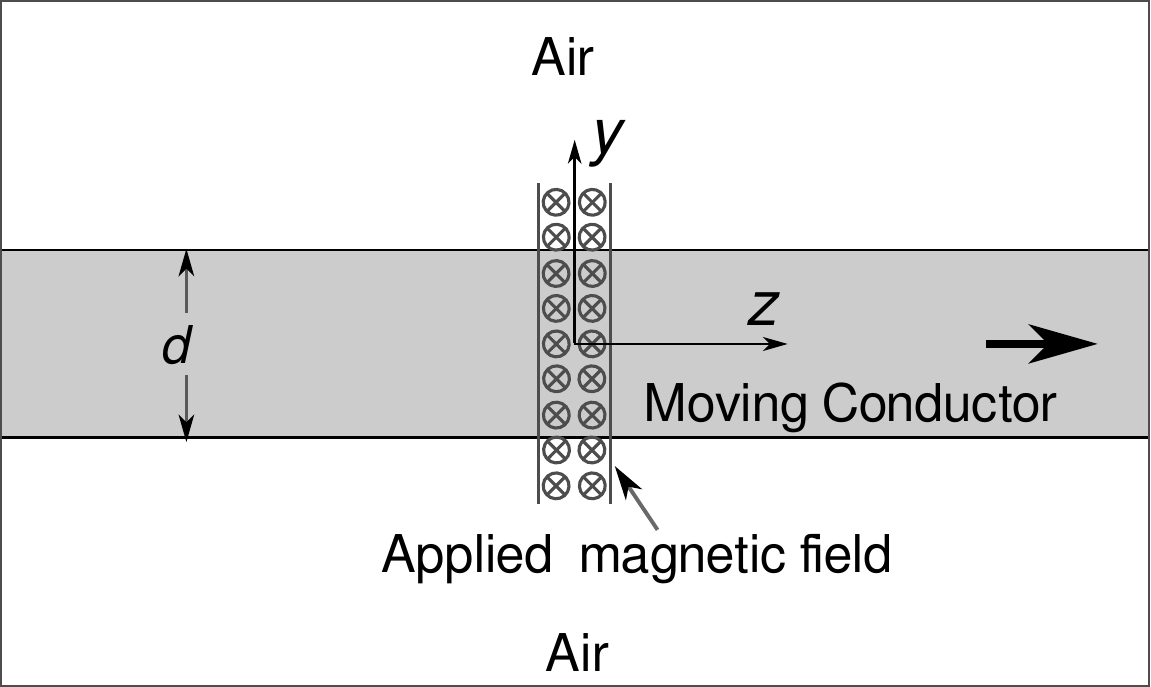}
		\caption{Schematic of the 2D problem.}
		\label{system2D}
\end{figure}

Stability analysis of a complete moving conductor problem is very difficult to handle, mainly due to the presence of multiple materials and the structure of the simulation domain. Therefore, a simplified moving conductor problem will be considered here, following a previous work\cite{su1}. A slightly modified version of the 2D moving conductor problem  used in \cite{su1} is shown in Fig.\ref{system2D}. In this, a conducting slab of thickness $d$ is moving along the $z$-axis with velocity $u_z$, under the influence of magnetic field $B_x$ directed along the $x$-axis. The conductivity and permittivity for the conductor are denoted as $\sigma$ and $\mu$ respectively. 

\begin{figure}
		\centering
		\mbox{\subfloat[]{\label{f:1dpr_2dshp} 
		\includegraphics[scale=0.7]{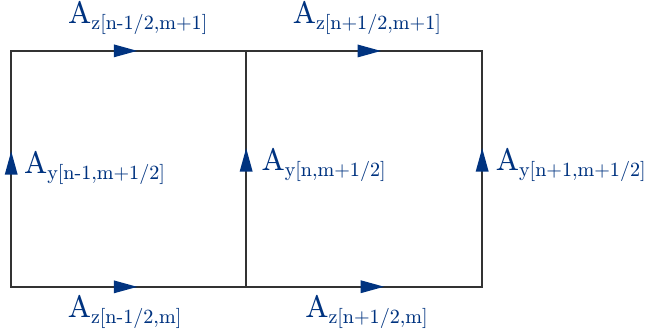}}}\\
		\mbox{\subfloat[]{\label{f:1dpr_1dshp} 
		\includegraphics[scale=0.7]{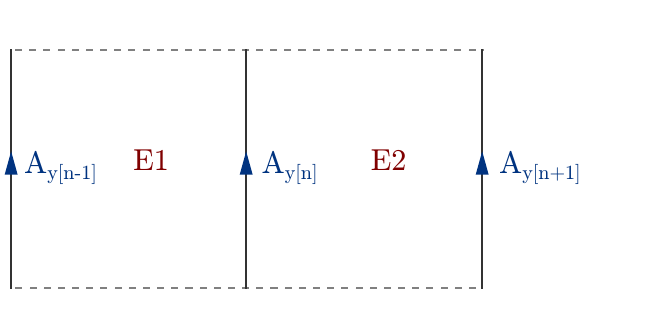}}}

		\caption{Representation of $\bf{A}$ with edge elements in $zy$-plane (a) 2D (b) reduced 1D.}
		\label{f:1dpr_shp}
\end{figure}

For this 2D problem, the vector potential has components of $A_y$ and $A_z$. The same has been depicted in Fig.\ref{f:1dpr_2dshp}. In Fig.\ref{f:1dpr_2dshp}, the finite element discretisation of the 2D problem using edge elements is shown for one $y$-edge ([n,m+1/2]). For the sake of mathematical analysis, a simplified version with equal discretisation along the $z$ and $y$ axis is chosen. The edge variables are subscripted with $n$ and $m$, where $n$ denotes the progression along the $z$ axis and $m$ denotes the progression along the $y$ axis. It can be noted that, in addition to the integer progression ($n-1$, $n$, $n+1$), a factor of $1/2$ is present to denote the edge variable that is constant for the edge.

Now, let us consider the limiting case of $d\to\infty$ as described in \cite{su1}. Here, due to the symmetry along $y$-axis, the variations with respect to $y$-axis vanishes, resulting in a problem which is independent of $\phi$ and $A_z$. This situation is depicted in Fig.\ref{f:1dpr_1dshp}, wherein only $A_y$ has variations along the $z$-axis as is the case for $y$ directed edges which have the natural variation along the  $z$-axis (perpendicular axis). The corresponding finite element formulation using the edge elements can be written as, 
\begin{equation}\label{1DFEMformulation}
\begin{split}
\int_\Omega \dfrac{d{M_y^l}}{dz} \frac{ d{{A}_y}}{dz} ~d\Omega + \mu \sigma u_z \int_\Omega          M_y^l       \frac{ d{{A}_y}}{dz} ~d\Omega ~\dots \\= \mu \sigma u_z\int_\Omega       M_y^l          {B}_x
\end{split}
\end{equation}
{where, $M_y^l$ is the $y$-directed edge weight function} \cite{emfembookjin}.
Evaluating the above equation for the $n^{th}$ edge gives the following difference equation,
\begin{equation}\label{eq:1Ddiff}
\begin{split}
-(1+P_e)A_{y[n-1]} + 2A_{y[n]} - (1-P_e)A_{y[n+1]} ~\dots \\= \frac{1}{3}(B_{x[n-1]} + 4*B_{x[n]} + B_{x[n+1]})*P_e*\Delta z
\end{split}
\end{equation}
This is same as that of the nodal formulation as in \cite{su2}. The numerical instability due to the negative roots (\emph{poles}) of the difference equation can also be viewed with the help of $z$-transform \cite{su1,su2}. Moreover, the $z$-transform clearly shows the effect of \emph{zeros} arising from the source term as well.
Applying the $Z$-transform on (\ref{eq:1Ddiff}), 
\begin{equation}\label{1Dexpr}
\frac{A_y}{B_x} = \frac{(Z+0.27)(Z+3.73)P_e\Delta z}{3(-1+P_e)(Z-1)(Z-\dfrac{-1-P_e}{-1+P_e})}
\end{equation}
For $P_e >>$ 1, the above equation (\ref{1Dexpr}) reduces to,
\begin{equation}\label{1DexprhighPe}
\frac{A_y}{B_x} \simeq \frac{\Delta z}{3} \frac{(Z+0.27)(Z+3.73)}{(Z-1)(Z+1)}
\end{equation}	
Now, the pole at -1 of equation (\ref{1DexprhighPe}) leads to the oscillation in the solution. 

Here, the relation between $A_y$ and $B_x$ in (\ref{1DexprhighPe}) is found to be same as that of in \cite{su2} with the nodal formulation. The oscillation appearing in the solution in \cite{su2} is successfully mitigated by applying the input field $B_x$ as the averaged nodal flux densities for the element. Here, the same approach is extended for the edge formulation, with elemental input fields, $B_{xE1}$, $B_{xE2}$.
The $B_{xE1}$ and $B_{xE2}$ are the element-averaged input magnetic fields for the elements spanning $[n-1], [n]$ and $[n],[n+1]$ respectively. The $B_{xE1}$ and $B_{xE2}$ are defined as \cite{su2},

\[ B_{xE1} = ( B_{x[n-1]} + B_{x[n]} ) /2 \]
\[ B_{xE2} = ( B_{x[n]}   + B_{x[n+1]})/2\]
For these modified inputs, the difference equation for (\ref{1DFEMformulation}) becomes, 
\begin{equation}
\begin{split}
-(1+P_e)A_{y[n-1]} + 2A_{y[n]} - (1-P_e)A_{y[n+1]} ~\dots \\= \frac{1}{2}(B_{x[n-1]} + 2*B_{x[n]} + B_{x[n+1]})*P_e*\Delta z
\end{split}
\end{equation}
Applying the $Z$-transform,
\begin{equation}\label{1Dexpravgip}
\frac{A_y}{B_x} = \frac{(Z+1)^2 P_e\Delta z}{2(-1+P_e)(Z-1)(Z-\dfrac{-1-P_e}{-1+P_e})}
\end{equation}
For $P_e >>$ 1, the above equation (\ref{1Dexpravgip}) reduces to,
\begin{equation}\label{1DexpravgiphighPe}
\frac{A_y}{B_x} \simeq \frac{\Delta z}{2} \frac{(Z+1)}{(Z-1)}
\end{equation}	

While comparing (\ref{1DexpravgiphighPe}) with (\ref{1DexprhighPe}), one would readily recognize that for high Peclet number, the zero at -1 for the second case eventually cancels the oscillatory pole at -1, thus leading to a stable solution. This indicates that the elemental average of the input magnetic field can be used for the edge elements as well. 
However, further confirmation in 2D would be helpful, since the edge elements are generally designed for curl problems.
The stability analysis in 2D is dealt in the next subsection.

%-------------------------------------------------------------------------------
%     Part 5.2. 2D Analysis
%-------------------------------------------------------------------------------
\subsection{Analysis with the 2D version of the problem}

\begin{figure}
\includegraphics[scale=0.37]{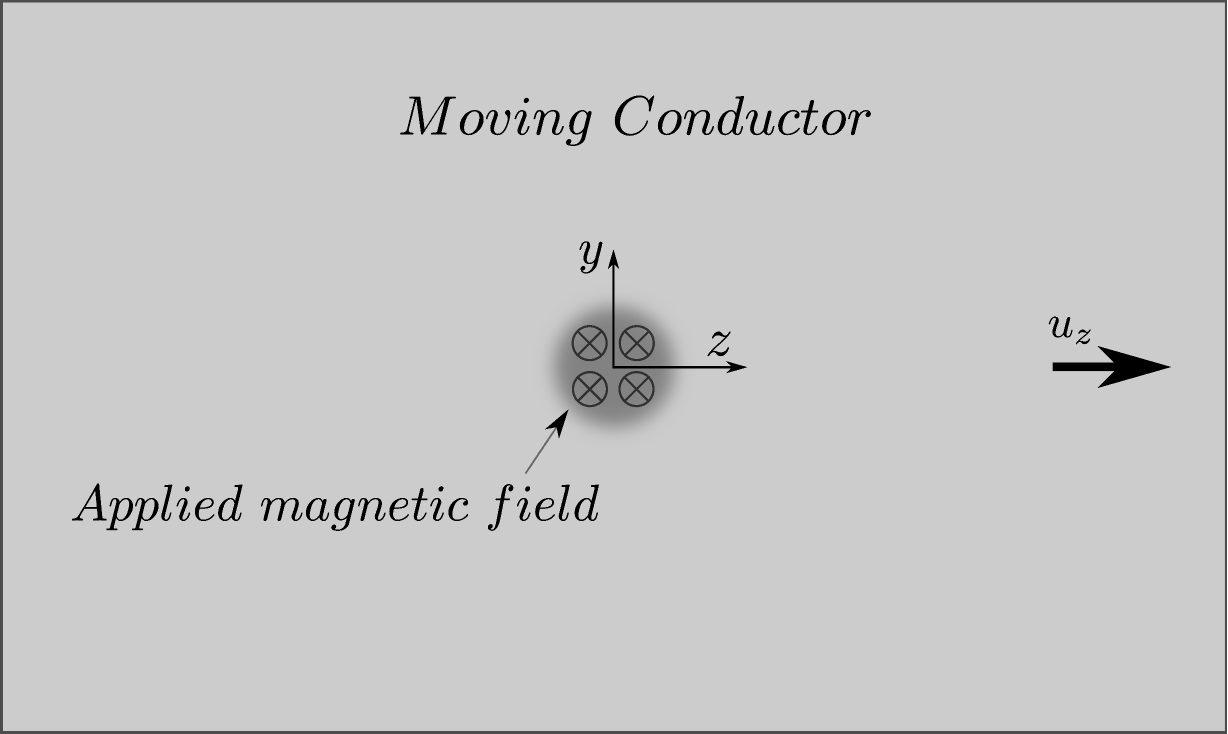}
\caption{Schematic of the simplified 2D problem for the stability analysis}
\label{f:2d_z}
\end{figure}

A simplified, 2D version of the moving conductor problem as shown in Fig. \ref{f:2d_z} is considered. The governing equations for the same are given below,
\begin{equation} \label{ge2dphi}
\nabla \cdot (\sigma \nabla \phi) - \nabla \cdot (\sigma {\bf u} \times (\nabla \times {\bf A})) = \nabla \cdot (\sigma u_z B_x {\bf \hat y})
\end{equation}
\begin{equation} \label{ge2da}
\sigma \nabla \phi - \nabla \cdot \dfrac{1}{\mu}(\nabla {\bf A}) - \sigma {\bf u} \times (\nabla \times {\bf A}) = \sigma u_z B_x {\bf \hat y}
\end{equation}

In the edge element formulation, the magnetic vector potential $\bf{A}$ is modeled with the edge vector shape functions $\bf{M}$
\cite{emfembookjin}. For this cartesian case, the $y$-component of the vector potential would be modeled with $y$-directed edge shape functions, and similarly, the $z$-component of the vector potential would be modeled with $z$-directed edge shape functions. The electric scalar potential $\phi$ is modeled with the nodal shape functions $N$. In the Galerkin finite element formulation, the weight functions are the shape functions themselves, and to mark the difference the weight functions are super-scripted with $l$. 

%In order to formulise for the 2D case, the $y$ and $z$ component of the magnetic vector potential $\bf{A}$ edge weight functions $\bf{M}$, i.e. $M_y^l$ and $M_z^l$ are considered. Here, $M_y^l$ is having a linear  variation along $z$-direction in the range $0$ to $1$, while being constant along $y$-direction. Similarly $M_z^l$ is having a linear  variation along $y$-direction in the range $0$ to $1$, while being constant along $z$-direction.

The Galerkin finite element formulation for the equations (\ref{ge2dphi}), (\ref{ge2da}) can be written as,
\begin{equation}\label{phiGE}
\begin{split}
\int_\Omega\dfrac{d{N^l}}{dz}\dfrac{d{\phi}}{dz}~d\Omega +
\int_\Omega\dfrac{d{N^l}}{dy}\dfrac{d{\phi}}{dy}~d\Omega + 
u_z \int_\Omega\dfrac{d{N^l}}{dy}\dfrac{d{A_y}}{dz}~d\Omega ~\dots \\
- u_z \int_\Omega\dfrac{d{N^l}}{dy}\dfrac{d{A_z}}{dy}~d\Omega
= u_z\int_\Omega\dfrac{d{N^l}}{dy} {B}_x
\end{split}
\end{equation}

\begin{equation}\label{yGE}
\begin{split}
\mu \sigma \int_\Omega M_y^l \dfrac{d{\phi}}{dy}~d\Omega + 
\int_\Omega \dfrac{d{M_y^l}}{dz}   \frac{ d{{A}_y}}{dz} ~d\Omega +
\int_\Omega \dfrac{d{M_y^l}}{dy}   \frac{ d{{A}_y}}{dy} ~d\Omega ~\dots \\+
  \mu \sigma u_z \int_\Omega M_y^l \frac{ d{{A}_y}}{dz} ~d\Omega  -
  \mu \sigma u_z \int_\Omega M_y^l \frac{ d{{A}_z}}{dy} ~d\Omega ~\dots \\
= \mu \sigma u_z \int_\Omega M_y^l {B}_x
\end{split}
\end{equation}

\begin{equation}\label{zGE}
 \mu\sigma \int_\Omega M_z^l \dfrac{d{\phi}}{dz}~d\Omega +
\int_\Omega \dfrac{d{M_z^l}}{dz}   \frac{ d{{A}_z}}{dz} ~d\Omega +
\int_\Omega \dfrac{d{M_z^l}}{dy}   \frac{ d{{A}_z}}{dy} ~d\Omega = 0
\end{equation}

It can be noted that, the equations (\ref{yGE}), (\ref{zGE}) are the weighted residual formulation of the equation (\ref{ge2da}); the former arise from the $y$-directed edge weight function ${\bf M}^l = M_y^l$ and the latter arise from the $z$-directed edge weight function ${\bf M}^l = M_z^l$. In addition to this, for the rest of the analysis, the element lengths along the $y$ and $z$-directions are assumed to be equal i.e. $\Delta y = \Delta z$ \cite{su2}.

\begin{figure*}
	\centering
	\mbox{\subfloat[]{\label{edge2d} 
    \includegraphics[scale=0.7]{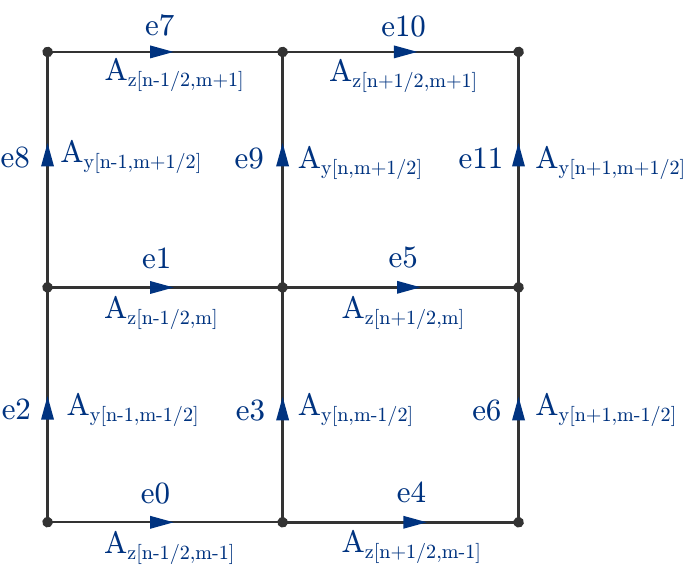}}} \hspace{5mm}
	\mbox{\subfloat[]{\label{node2d} 
	\includegraphics[scale=0.7]{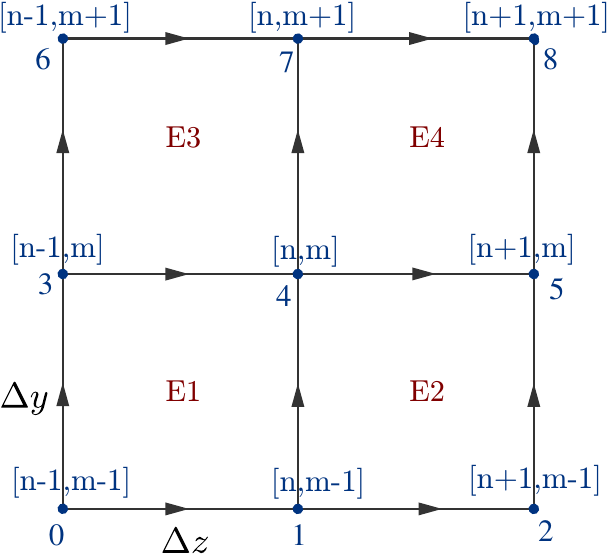}}}

    \caption{Representation of edge elements in 2D $zy$-plane (a) $z$ and $y$ directed edge vectors and their corresponding edge numbers (b) node and element numbering for the same set of edge elements.}
    \label{2dA}
\end{figure*}

The difference form of the finite element equation (\ref{phiGE}) for node-4 (see Fig. \ref{2dA}) can be written as,   
\begin{equation}\label{NL}
\begin{split}
\dfrac{1}{6}(2(\phi_{[0]} + \phi_{[1]} + \phi_{[2]} +\phi_{[3]} - 8\phi_{[4]} + \phi_{[5]} + \phi_{[6]}~\dots \\
 + \phi_{[7]} + \phi_{[8]}) 
+ 3u_z (A_{[e2]} - A_{[e0]} - A_{[e4]} - A_{[e6]}~\dots \\
 + A_{[e10]} - A_{[e8]} + A_{[e11]} - A_{[e7]}))~\dots \\
= \dfrac{u_z\Delta y }{12}(B_{x[8]} + B_{x[6]} + 4B_{x[7]} - B_{x[0]}~\dots \\ - 4B_{x[1]} - B_{x[2]})
\end{split}
\end{equation}

In the above equation, it may be noted that the scalar variable $\phi$ is at the nodes, whereas, the vector potential $A_{[e]}$ is present at the edges. The stability analysis requires the variables to be at the nodes \cite{su2, 2dhb, 2dsb}. Therefore, it is necessary to represent the vector potential at the nodes. For this, the node equivalent of the vector potential is obtained by averaging the vector potential at the edges. Each node has two $y$-directed edges (vertical edges); by taking the average of the two, the $A_{y[node]}$ is obtained. Similarly, each node has two $z$-directed edges (horizontal edges); by taking the average of the two, the $A_{z[node]}$ is obtained. With these, the nodal form of the (\ref{NL}) can be written as, 
\begin{equation}\label{NLnode}
\begin{split}
\dfrac{1}{6}(2(\phi_{[0]} + \phi_{[1]} + \phi_{[2]} +\phi_{[3]} - 8\phi_{[4]} + \phi_{[5]} + \phi_{[6]}~\dots \\
 + \phi_{[7]} + \phi_{[8]}) 
+ \dfrac{3u_z}{2} (A_{y[0]} + A_{y[8]} - A_{y[2]} - A_{y[6]})~\dots \\
 + 3u_z(-A_{z[1]} + A_{z[7]}))~\dots \\
= \dfrac{u_z\Delta y}{12}(B_{x[8]} + B_{x[6]} + 4B_{x[7]} - B_{x[0]}~\dots \\ - 4B_{x[1]} - B_{x[2]})
\end{split}
\end{equation}

For the vector equation, the 2D representation (see Fig. \ref{edge2d}) has two edges for each direction. That is, the shape function corresponds to edges $e3,~e9$, each form a finite element equation. Similarly, the shape function corresponds to edges $e1,~e5$, each form a finite element equation. Thus, there are 4 finite element equations associated with 4-edges. In contrast, the node element form only one equation that corresponds to node-4. The edges $e1,~e5$ are directed along the $z$ direction; so their Galerkin weighted residual formulation corresponds to the $z$ component of the vector equation (\ref{ge2da}). Then, the difference form of the finite element equation (\ref{zGE}) for $e1$ can be written as,   
\begin{equation}\label{ML1}
\begin{split}
-\dfrac{1}{6}((\phi_{[0]} - \phi_{[1]} + 4\phi_{[3]} - 4\phi_{[4]} + \phi_{[6]} - \phi_{[7]})\Delta y^2 \mu \sigma ~\dots \\
+ 6(A_{[e0]} - 2A_{[e1]} + A_{[e7]})\Delta z)/(\Delta y^2 \Delta z \mu) = 0 
\end{split}
\end{equation}
Similarly, the difference form of the finite element equation (\ref{zGE}) for $e5$ can be written as,   
\begin{equation}\label{ML5}
\begin{split}
-\dfrac{1}{6}((\phi_{[1]} - \phi_{[2]} + 4\phi_{[4]} - 4\phi_{[5]} + \phi_{[7]} - \phi_{[8]})\Delta y^2 \mu \sigma ~\dots \\
+ 6(A_{[e10]} - 2A_{[e5]} + A_{[e4]})\Delta z)/(\Delta y^2 \Delta z \mu) = 0 
\end{split}
\end{equation}
Upon taking the average of equations (\ref{ML1}) and (\ref{ML5}), one gets the averaged difference form of the finite element equation (\ref{zGE}) along $z$-direction. 
\begin{equation}\label{ML1ML5}
\begin{split}
-\dfrac{1}{6}((\phi_{[0]} - \phi_{[2]} + 4\phi_{[3]} - 4\phi_{[5]} + \phi_{[6]} - \phi_{[8]})\Delta y^2 \mu \sigma ~\dots \\
+ 6(A_{[e0]} -2A_{[e1]} + A_{[e4]} - 2A_{[e5]} + A_{[e7]} ~\dots \\
+ A_{[e10]})\Delta z)/(\Delta y^2 \Delta z \mu) = 0 
\end{split}
\end{equation}
The above equation (\ref{ML1ML5}) is in terms of the edge vectors $A_{[e]}$. By following the similar procedure that applies for (\ref{NLnode}), the edge variables are represented with their respective node equivalents,
\begin{equation}\label{ML1ML5node}
\begin{split}
-\dfrac{1}{6}((\phi_{[0]} - \phi_{[2]} + 4\phi_{[3]} - 4\phi_{[5]} + \phi_{[6]} - \phi_{[8]})\Delta y^2 \mu \sigma ~\dots \\
+ 12(A_{z[1]} -2A_{z[4]} + A_{z[7]})\Delta z)/(\Delta y^2 \Delta z \mu) = 0 
\end{split}
\end{equation}

By following the same procedure for the $y$ directed edges; that is i) obtain the difference equations of the Galerkin formulation from $e3$ and $e9$ edges ii) average them iii) represent the edge variables with their node equivalents. After these 3-steps, the final difference equation for (\ref{yGE}) can be written as,

\begin{equation}\label{ML3ML9node}
\begin{split}
-\dfrac{1}{6}((\phi_{[0]} + 4\phi_{[1]} + \phi_{[2]} - \phi_{[6]} - 4\phi_{[7]} - \phi_{[8]})\Delta z \mu \sigma ~\dots \\
-6(-A_{y[3]} +  A_{y[5]} + A_{z[1]} -  A_{z[7]})\mu \sigma u_z \Delta z ~\dots \\
+ 12(A_{y[3]} -2A_{y[4]} + A_{y[5]}))/(\Delta z^2 \mu)~\dots \\
= \dfrac{1}{12}(B_{x[0]} + 4B_{x[1]} + B_{x[2]} + 2B_{x[3]} +  8B_{x[4]} ~\dots \\
+ 2B_{x[5]} + B_{x[6]} + 4B_{x[7]} + B_{x[8]})\sigma u_z 
\end{split}
\end{equation}

The equations (\ref{NLnode}), (\ref{ML1ML5node}), (\ref{ML3ML9node}) forms a system of 3 equations and 3 variables $\phi, A_y, A_z$. The stability analysis can be now performed by taking 2D $Z$-transform of these equations \cite{2dhb, dcbook1, su2}. In this, the $z$-direction is represented by the transformation
variable $Z_n$ and the $y$-direction is represented by the transformation
variable $Z_m$ (please refer to Fig. \ref{node2d}). After taking 2D $Z$- transform, the equations (\ref{NLnode}), (\ref{ML3ML9node}), (\ref{ML1ML5node}) take the following form,
\begin{equation}\label{NLZ}
\begin{split}
\dfrac{1}{3}[S1]\phi + \dfrac{u_z}{4}[S2]A_y - \dfrac{u_z}{2}[S3]A_z
= \dfrac{u_z \Delta z}{12}[Q1]B_x
\end{split}
\end{equation}
\begin{equation}\label{MLyZ}
\begin{split}
\dfrac{Pe}{6u_z}[Q1]\phi + (-[S1]+Pe[Q2])A_y + Pe[S3]A_z ~\dots \\
=\dfrac{Pe \Delta z}{12}[M1]B_x
\end{split}
\end{equation}
\begin{equation}\label{MLzZ}
\begin{split}
\dfrac{Pe}{6u_z}[Q2']\phi -[S1']A_z = 0
\end{split}
\end{equation}
where, 
\begin{equation}\label{Q1}
[Q1] = Z_n ^2 Z_m ^2 + 4 Z_n Z_m ^2 + Z_m ^2 - Z_n ^2 - 4 Z_n - 1 
\end{equation}
\begin{equation}\label{Q2}
[Q2] = Z_n ^2 Z_m - Z_m
\end{equation}
\begin{equation}\label{Q2'}
[Q2'] = Z_n ^2 Z_m ^2 - Z_m ^2 + 4 Z_n ^2 Z_m - 4 Z_m + Z_n ^2 - 1 
\end{equation}
\begin{equation}\label{S1}
\begin{split}
[S1] = 1 + Z_n + Z_n ^2 + Z_m + Z_n ^2 Z_m + Z_m ^2 ~\dots \\
+ Z_n Z_m ^2 + Z_n ^2 Z_m ^2 - 8 Z_n Z_m 
\end{split}
\end{equation}
\begin{equation}\label{S2}
[S2] = 1 + Z_n ^2 Z_m ^2 - Z_n ^2 - Z_m ^2
\end{equation}
\begin{equation}\label{S3}
[S3] = Z_n - Z_n Z_m ^2
\end{equation}
\begin{equation}\label{S1'}
[S1'] = Z_n - 2 Z_n Z_m + Z_n Z_m ^2
\end{equation}
\begin{equation}\label{M1}
\begin{split}
[M1] = 1 + 4 Z_n + Z_n ^2 + 2 Z_m + 8 Z_n Z_m~\dots \\  
+ 2 Z_n ^2 Z_m + Z_m ^2 + 4 Z_n Z_m ^2 + Z_n ^2 Z_m ^2
\end{split}
\end{equation}

The equations (\ref{NLZ}), (\ref{MLyZ}) and (\ref{MLzZ}) form a system of 3 equations with 3 variables $\phi, A_y, A_z$. By following a similar procedure
as that of \cite{su2}, the final equation is obtained for $Pe >> 1$ as,
%\begin{equation}\label{AyBxgen1}
%\begin{split}
%((2[S1][S1'] - \dfrac{Pe}{2}[Q2'][S3])(-[S1] + Pe[Q2]) ~\dots \\
%- ([Q1][S1'] + Pe[Q2'][S3])\dfrac{Pe}{4}[S2])A_y        ~\dots \\
%= ((2[S1][S1'] - \dfrac{Pe}{2}[Q2'][S3])(\dfrac{Pe\Delta z}{12} [M1]) ~\dots \\
%- ([Q1][S1'] + Pe[S3])(\dfrac{Pe\Delta z}{12} [Q1])) B_x
%\end{split}
%\end{equation}
%
%For $Pe>>1$ the above equation (\ref{AyBxgen1}) takes the following form,
%\begin{equation}\label{AyBxgen}
%\begin{split}
%((-\dfrac{Pe}{2}[Q2'][S3])Pe[Q2] - Pe[Q2'][S3]\dfrac{Pe}{4}[S2])A_y ~\dots \\
%\approx (-\dfrac{Pe}{2}[Q2'][S3]\dfrac{Pe\Delta z}{12}[M1] - Pe[S3]\dfrac{Pe\Delta z}{12}[Q1])B_x
%\end{split}
%\end{equation}

\begin{equation}\label{AyBxgen}
\begin{split}
\left([Q2][S3]) +  \dfrac{[S3][S2]}{2} \right) A_y 
\approx \dfrac{\Delta z}{12} (~ [M1][S3]  ~\dots \\ 
                                + 2 [Q1][S3] ~)B_x
\end{split}
\end{equation}
Upon simplification, the transfer function takes the following form,
\begin{equation}\label{AyBxTF}
\dfrac{A_y}{B_x} \approx \dfrac{\Delta z}{6} \dfrac{(Z_n^2 + 4 Z_n + 1)f1(Z_m)}{(1-Z_n^2)f2(Z_m)}
\end{equation}
where,
\[ f1(Z_m) = 1 - 2 Z_m - 3 Z_m^2 \]
\[ f2(Z_m) = (2 Z_m - (1-Z_m^2)^2) \]

The transfer function (\ref{AyBxTF}) has \emph{poles/roots} located at
 $-1, +1$, among which the pole at `-1' causes the oscillatory behaviour in the solution. The \emph{zeros} are located at $Z_n=-0.27, -3.7$. The locations of the
zeros may be changed by modifying the representation of the input magnetic
field consistently. By following the 1D case (please refer to \ref{sec1D.A})
and the preceding work \cite{su2}, the input magnetic field, which is 
directed along $x$-edges (perpendicular to the plane) can be averaged over
each element and the resultant magnetic field can be used as the input field.
In other words, the following would be the input field for the element $E1$,
\begin{equation}\label{BxAvg2D}
\begin{split}
B_x = \dfrac{1}{4} (B_x[n-1,m-1] + B_x[n,m-1] \dots \\ 
                                 + B_x[n-1,m] + B_x[n,m])
\end{split}
\end{equation}
and similarly for other elements. With this modification, the transfer function for the proposed formulation is calculated to take the following form,
\begin{align}\label{AyBx_Bavg_TF}
\dfrac{A_y}{B_x} &\approx \dfrac{\Delta z}{4} \dfrac{(Z_n^2 + 2 Z_n + 1)f1(Z_m)}{(1-Z_n^2)f2(Z_m)} \nonumber \\
&\approx \dfrac{\Delta z}{4} \dfrac{(1 + Z_n)\bcancel{(1+Z_n)}f1(Z_m)}{(1 - Z_n)\bcancel{(1+Zn)}f2(Z_m)} 
\end{align}
The oscillatory pole at `$Z = -1$' is canceled by the zero introduced in the
numerator. Thus, the elemental averaging of the input magnetic field can be seen to possess the stability properties similar to that of 1D. In the next section,
numerical validation exercises are carried out in 2D and 3D.

%-------------------------------------------------------------------------------
% Numerical Simulation
%===============================================================================

\section{Numerical Validation}

%===============================================================================
%-------------------------------------------------------------------------------

\subsection{Simulation Results for 2D version of the problem}

%-------------------------------------------------------------------------------
\begin{figure}
		\centering
		\includegraphics[scale=0.2]{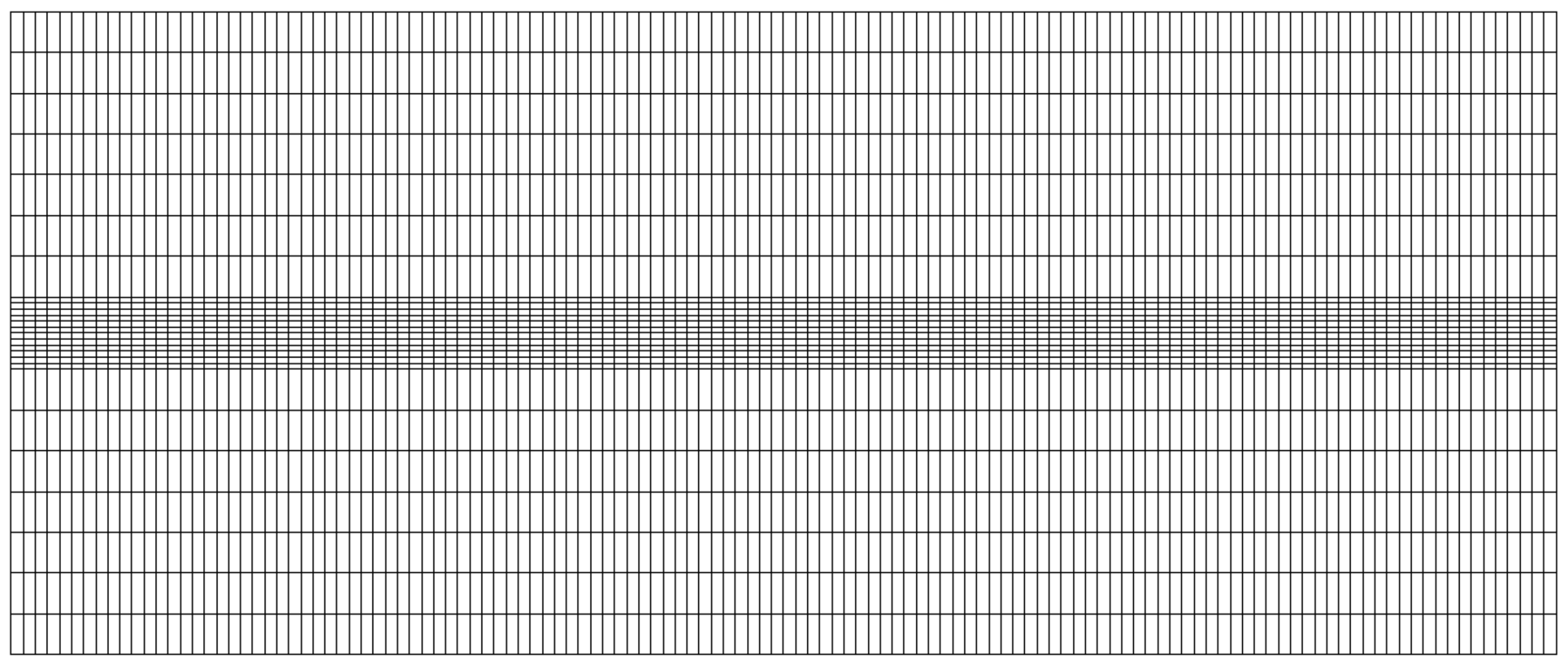}
		\caption{ Sample 2D finite element mesh }
		\label{2dmesh}
\end{figure}

\begin{figure*}
		\centering
		\mbox{\subfloat[]{\label{2D_MC_EdgeA} 
		\includegraphics[scale=0.6]{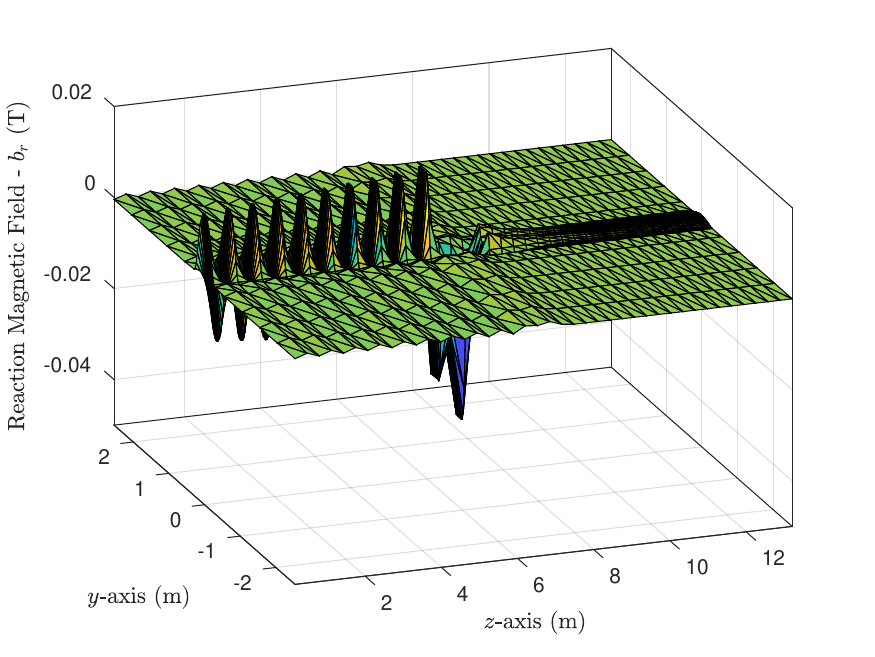}}}
		\mbox{\subfloat[]{\label{2D_MC_ElemA} 
		\includegraphics[scale=0.6]{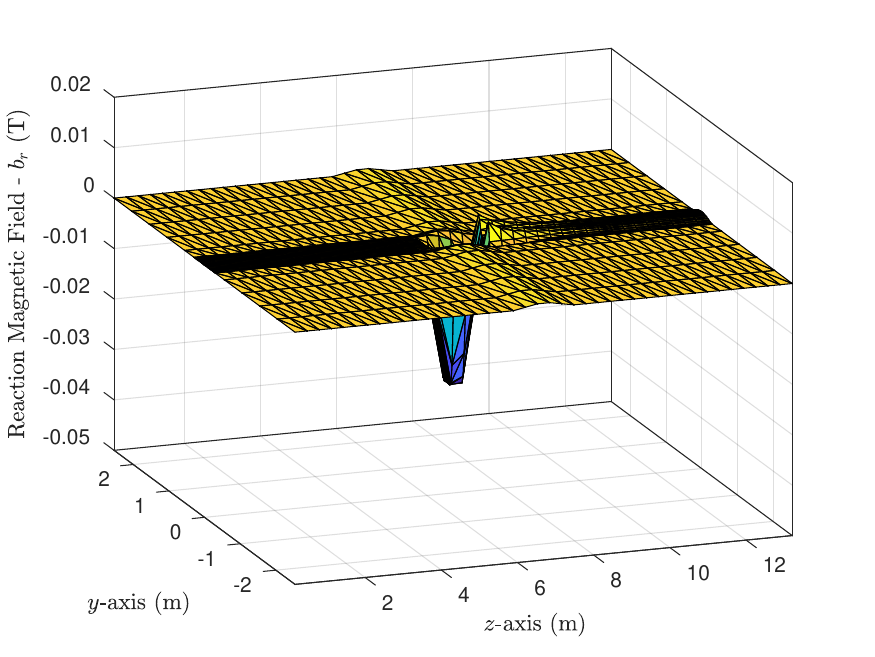}}}\\
		\caption{ Reaction magnetic field - $b_x$ from the 2D moving 
		conductor problem with $Pe=200$ 
		(a) Galerkin formulation
		(b) Proposed formulation}
		\label{2dMCvalid}
\end{figure*}

The 2D simulation involves the problem setup, shown in Fig. \ref{system2D}. 
The corresponding finite element mesh is shown in Fig. \ref{2dmesh}.
The physical parameters of the problem are as follows. The conductor has a width 
of $d = 0.5m$, its conductivity is 
$\sigma=7.2 \times 10^6 Sm^{-1}$ and its velocity is $u_z = 50 ms^{-1}$.
Simulations were carried out to test the stability of the proposed 
formulation; stable solutions are observed. 
A sample simulation result showing the reaction magnetic field 
$b_x = \partial A_z / \partial y - \partial A_y / \partial z$ with $Pe=200$ is displayed in Fig. \ref{2dMCvalid}.
The Fig. \ref{2D_MC_EdgeA} shows the 
$b_x$ obtained from the standard Galerkin formulation and the Fig. \ref{2D_MC_ElemA} shows the same from the proposed formulation.
It can be seen that the proposed formulation gives a stable solution without
any numerical oscillations.

In the nodal formulation, it was possible to obtain the analytical
expression for the peak errors due to the numerical oscillation \cite{su2}.
This is because, the nodal formulation can be reduced to 1D and the
resulting finite element equation in the difference form can be solved.
The following expressions of peak errors are derived in \cite{su2}.
\vspace{2mm}

Analytical error in the standard Galerkin formulation:
\begin{equation} \label{eq:errGA}
\hat{e}_{GA} = \left| \dfrac{(Pe^2 - 3)(Pe - 1)}{3(Pe + 1)^3} \right| \times 100 \%
\end{equation}

Analytical error in the stable nodal formulation of \cite{su2}:
\begin{equation} \label{eq:errGA}
\hat{e}_{SA} = \left| \dfrac{Pe - 1}{(Pe + 1)^3}  \right| \times 100\%
\end{equation}

\begin{figure}
		\centering
		\includegraphics[scale=0.62]{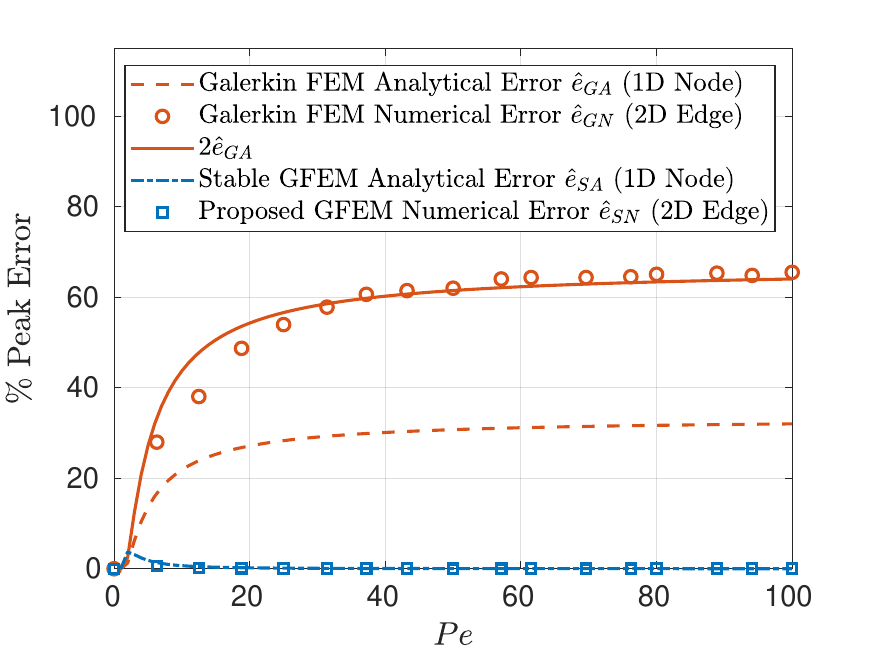}
		\caption{\% Peak error measured}
		\label{2dpkerr}
\end{figure}

Having these nodal errors as reference, the peak oscillation error with the 
edge elements are measured for the 2D problem. These measured values are 
plotted in Fig. \ref{2dpkerr}, where $\hat{e}_{GN}$ is the peak error measured
with the Galerkin formulation and $\hat{e}_{SN}$ is the peak error measured
with the Proposed formulation. It is observed that the peak error from the
2D-Galerkin-Edge is twice that of the 1D-Galerkin-Node. In other words, 
`$\hat{e}_{GN} \approx 2~ \hat{e}_{GA}$' and the same can observed from the
plot of `$2~ \hat{e}_{GA}$' in Fig. \ref{2dpkerr}. It can also be observed 
that the peak error measured with the proposed formulation ($\hat{e}_{SN}$)
 is negligible.

\newcolumntype{C}[1]{>{\centering\let\newline\\\arraybackslash\hspace{0pt}}m{#1}}
\begin{table*} [tbp]%[tbp]
\sffamily
%\begin{minipage}{\textwidth}
\centering
\caption{Measured values of Error in the first derivative for Galerkin scheme and the proposed scheme}
\label{Tab:1}
{\renewcommand{\arraystretch}{1.6}
\begin{tabular}{@{}|C{0.1\textwidth}|C{0.1\textwidth}|C{0.1\textwidth}|C{0.1\textwidth}|C{0.1\textwidth}|C{0.1\textwidth}|C{0.1\textwidth}|C{0.1\textwidth}|}
\hline 
\multicolumn{1}{|C{0.1\textwidth}}{} \vline & \multicolumn{1}{C{0.1\textwidth}}{}\vline & \multicolumn{6}{C{0.6\textwidth}}{\hspace{-60pt} Error measured in $b_x = \nabla \times {\bf A} \cdot \hat{x}$, normalised w.r.t. applied magnetic field $B_{ax}$} \vline \\
\cline{3-8}
\multicolumn{1}{|C{0.1\textwidth}}{Number of Elements} \vline & \multicolumn{1}{C{0.1\textwidth}}{$P_e$}\vline & \multicolumn{3}{C{0.3\textwidth}}{Galerkin Formulation} \vline & \multicolumn{3}{C{0.3\textwidth}}{Proposed Formulation} \vline \\

\cline{3-8}
\multicolumn{1}{|C{0.1\textwidth}}{} \vline & \multicolumn{1}{C{0.1\textwidth}}{}\vline & \multicolumn{1}{C{0.1\textwidth}}{L2 Error} \vline & \multicolumn{1}{C{0.1\textwidth}}{Absolute Error} \vline & \multicolumn{1}{C{0.1\textwidth}}{Expt. Order of Convergence} \vline& \multicolumn{1}{C{0.1\textwidth}}{L2 Error} \vline & \multicolumn{1}{C{0.1\textwidth}}{Absolute Error} \vline & \multicolumn{1}{C{0.1\textwidth}}{Expt. Order of Convergence} \vline\\
 
\hline
 640   & 100 & 6.035e-04 & 1.213e-01 &   -     & 1.765e-04 & 1.594e-02 &   -\\

 \hline
 2560  & 50 & 4.634e-04 & 8.741e-02 & 0.47 & 1.210e-04 & 8.650e-03 & 0.88\\

 \hline
 10240 & 25  & 2.708e-04 & 4.327e-02 & 1.01 & 7.919e-05 & 3.988e-03 & 1.12\\

 \hline
 40960 & 12.5  & 1.218e-04 & 1.220e-02 & 1.83 & 4.558e-05 & 1.509e-03 & 1.40\\

 \hline
\end{tabular}}
%\end{minipage}
\end{table*} 

In table \ref{Tab:1}, the average and the rms error measured for the 
first derivative (reaction magnetic field $b_x$) is presented. The
errors measured for the Galerkin scheme and the proposed scheme, for 
different amount of resolution. The average error, as well as, the rms (L2) error are observed to fall with the increasing resolution. It may also be noted
that the errors measured from the proposed scheme are an order of magnitude
smaller than the errors measured from the Galerkin scheme. With the
increasing resolution, the proposed formulation also produces the expected
convergence rate. Thus, the proposed formulation gives stable, as well as,
accurate results. In the next subsection, further testing is carried out
in 3D with the `Testing Electromagnetic Analysis Methods' (TEAM) 
problem No. 9 \cite{team9}.

%-------------------------------------------------------------------------------

\subsection{Validation with 3D TEAM-9 problem}

%-------------------------------------------------------------------------------

\begin{figure*}
		\centering
		\mbox{\subfloat[]{\label{team9a_sch} 
		\includegraphics[scale=0.37]{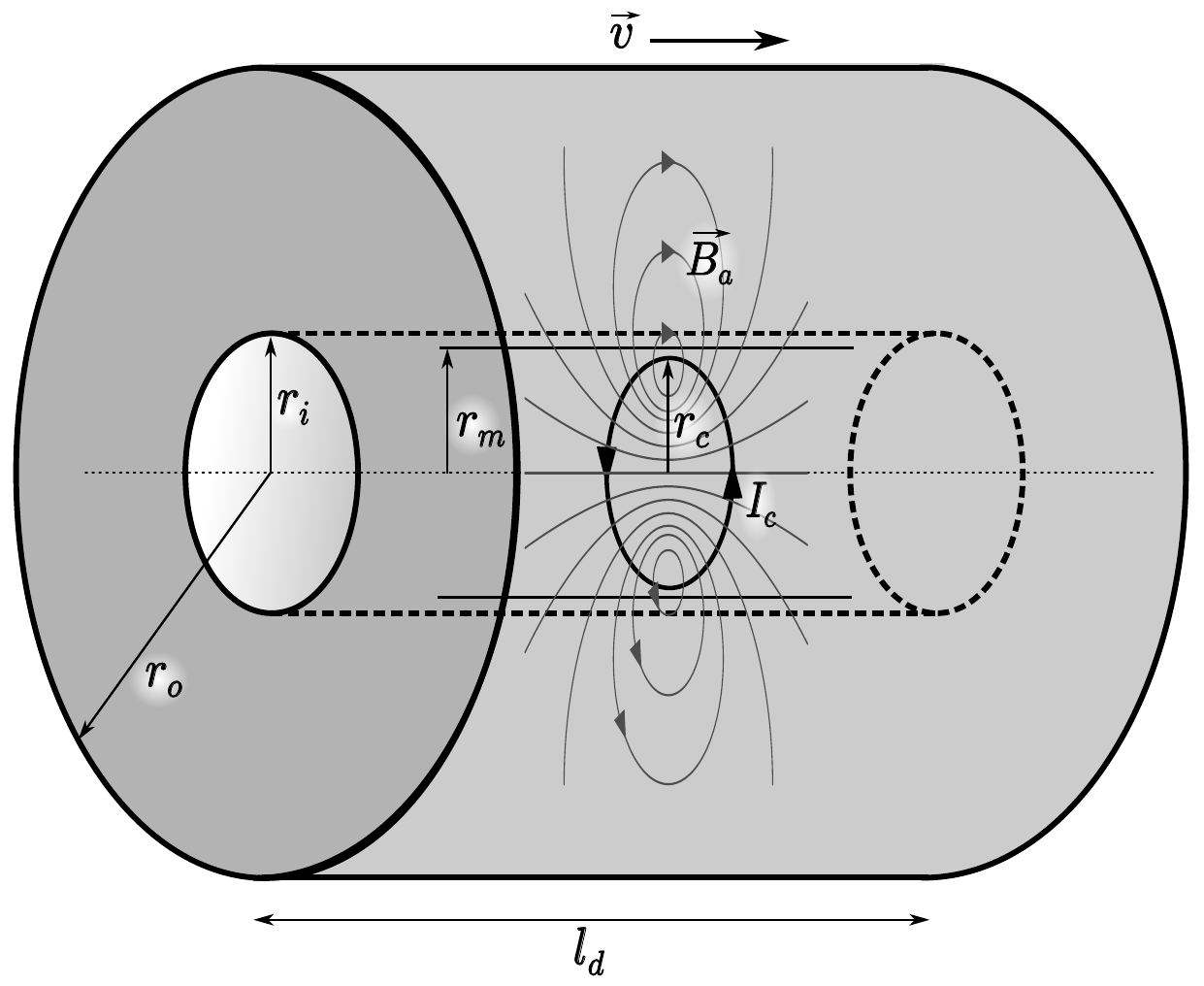}}}
		\mbox{\subfloat[]{\label{team9a_3dmesh} 
		\includegraphics[scale=0.24]{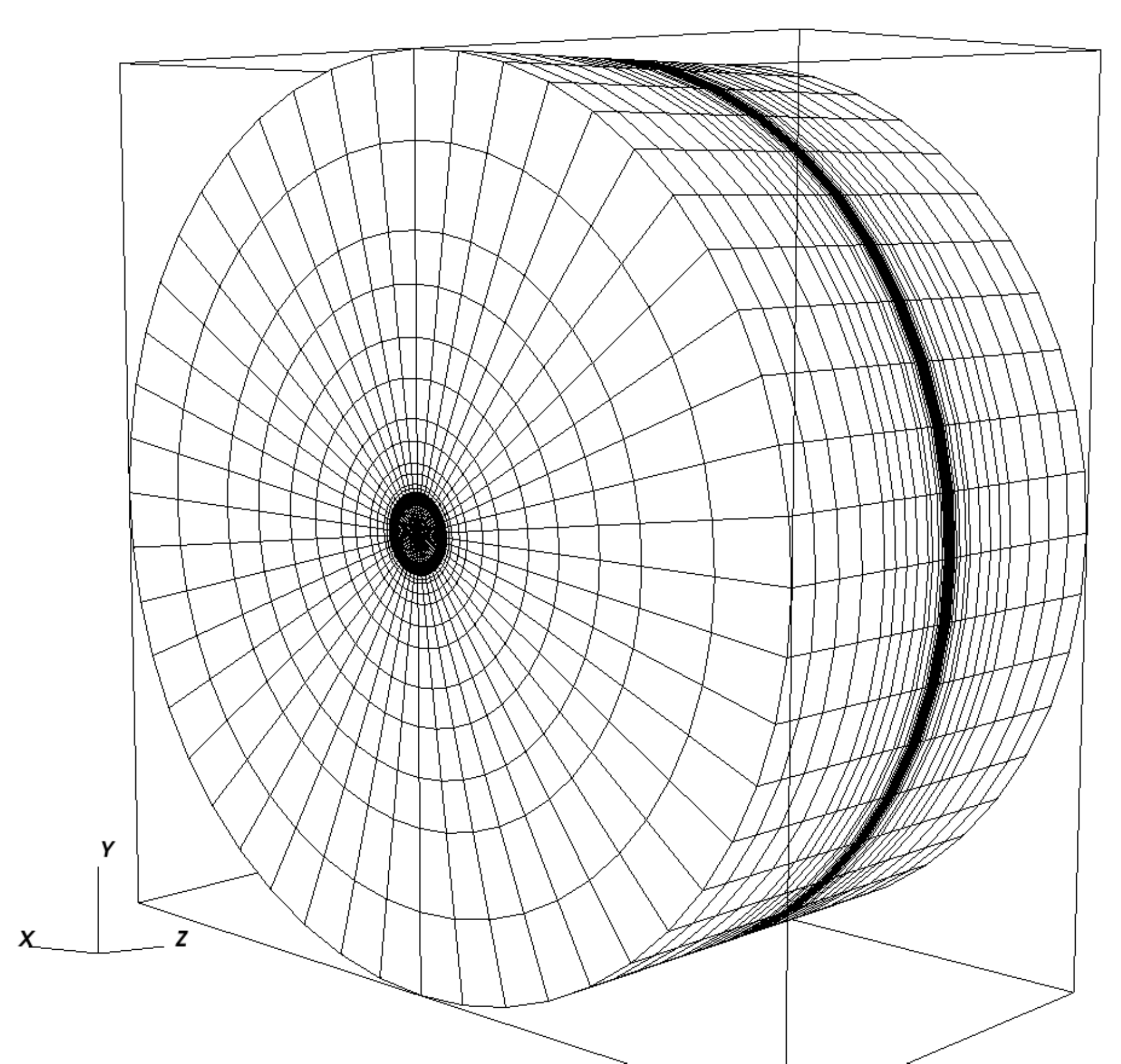}}}\\
		\mbox{\subfloat[]{\label{Edge_A} 
		\includegraphics[scale=0.6]{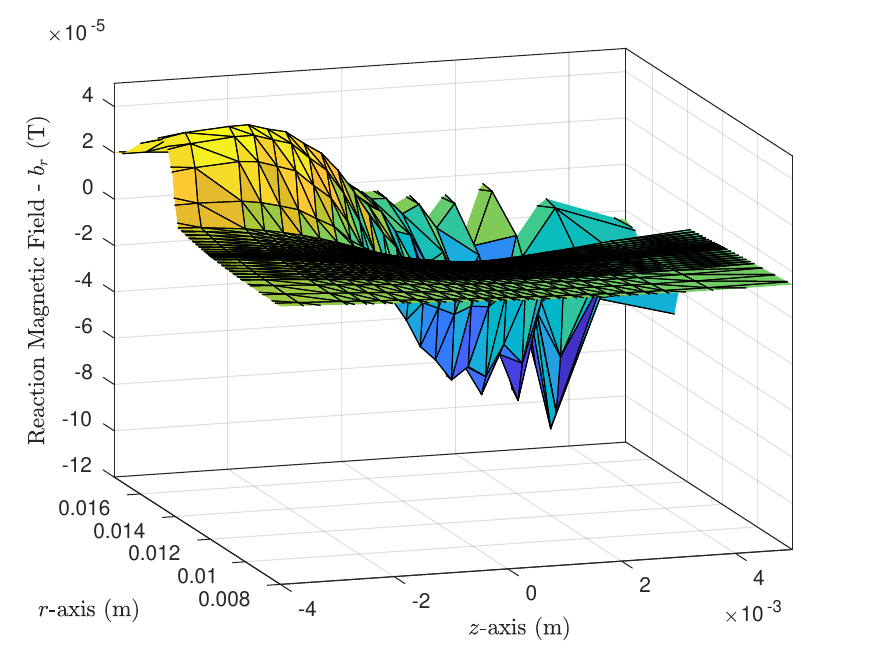}}}
		\mbox{\subfloat[]{\label{Elem_A} 
		\includegraphics[scale=0.6]{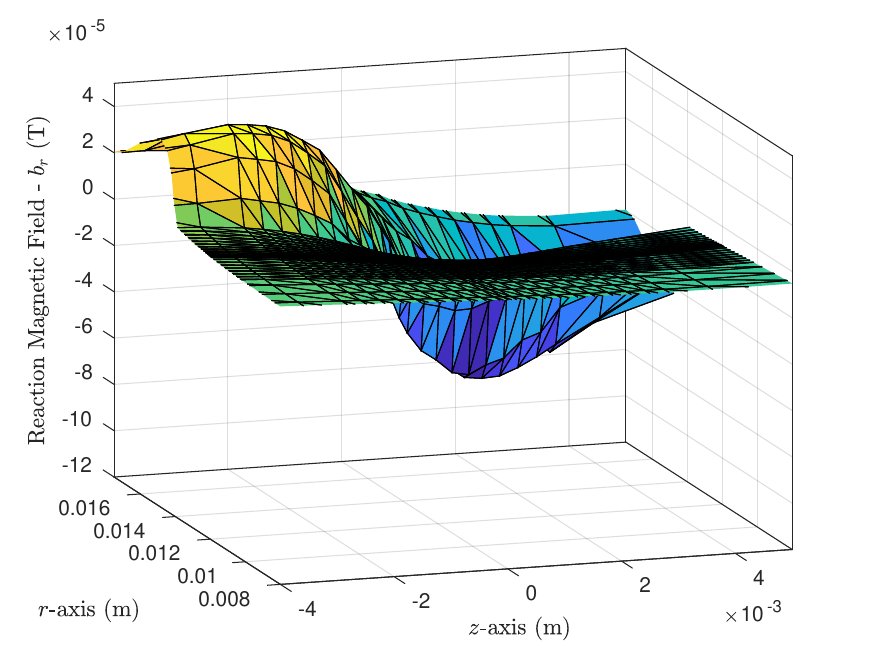}}}
		\caption{Description of the TEAM 9a problem and sample results
		(a) Schematic representation of the TEAM 9a moving conductor
		problem. 
		(b) Finite element mesh employed 
		(c) Galerkin scheme - reaction magnetic field - $b_r$ for $u_z$ = 100 $ms^{-1}$ and $\mu_r$ = 50 
		(d) Proposed scheme - reaction magnetic field - $b_r$ for $u_z$ = 100 $ms^{-1}$ and $\mu_r$ = 50  in the cross section along $rz$ plane for 
$\theta \approx 0^o$.}
		\label{3dteam9valid}
\end{figure*}

A schematic diagram of the TEAM-9 problem is shown in Fig. \ref{team9a_sch}.
The problem has an infinite ferromagnetic material with the conductivity of $\sigma=5\times 10^6~Sm^{-1}$. The relative magnetic permeability of the 
material is taken as $\mu_r=1$ and $\mu_r=50$. The ferromagnetic material has 
a cylindrical bore with the radius of $r_i = 14 \times 10^{-3} m$. 
A concentric current-carrying loop with $1A$ of current and a radius of $r_c = 12 \times 10^{-3} m$ is moving at an uniform velocity inside the bore.
For the analysis, the case with the largest velocity ($v=100ms^{-1}$) is
chosen. In order to accurately model the current loop, the finite element mesh
close to the current loop is dense; away from the current loop, the mesh 
becomes progressively coarser. Due to this variation, the resulting value of
the Peclet number varies from $5$ to $200$. The finite element mesh employed is shown in Fig. \ref{team9a_3dmesh}. In the proposed formulation, for each element, the $x, y, z$ components of the applied magnetic field $B_x, B_y, B_z$ are represented as,
\begin{equation}\label{bxavg3d}
B_x = \dfrac{1}{\sum_{e=1}^n |v^e_x| } \sum_{e=1}^{n} B^e  v^e_x 
\end{equation}
\begin{equation}\label{byavg3d}
B_y = \dfrac{1}{\sum_{e=1}^n |v^e_y| } \sum_{e=1}^{n} B^e  v^e_y 
\end{equation}
\begin{equation}\label{bzavg3d}
B_z = \dfrac{1}{\sum_{e=1}^n |v^e_z| } \sum_{e=1}^{n} B^e  v^e_z 
\end{equation}
where $n$ is the number of edges for each element, $B^e$ is the applied magnetic field corresponding to each edge, and $v_x^e, v_y^e, v_z^e$ form the unit vector ${\bf v}_e = \{v_x^e, v_y^e, v_z^e\}$ of the edge `$e$'. The elemental applied magentic field vector ${\bf B} = \{B_x, B_y, B_z \}$ is the source field in the proposed formulation. The standard Galerkin formulation would have the applied magnetic field at each Gauss-integration point ${\bf B}_g$ as,
\[ {\bf B}_g = \sum_{e=1}^n B^e {\bf M}_g^e \]
where ${\bf M}^e_g$ is the value of edge shape function vector at a Gauss-integration point.
The simulated, reaction magnetic field along the 
$rz$-plane is plotted in Fig. \ref{Edge_A} and Fig. \ref{Elem_A} for the
Galerkin scheme and the proposed formulation, respectively. It may be readily
noted that the $b_r$ from the proposed formulation is stable as expected. 

\begin{figure}
		\centering
		\mbox{\subfloat[]{\label{t9mur1} 
		\includegraphics[scale=0.62]{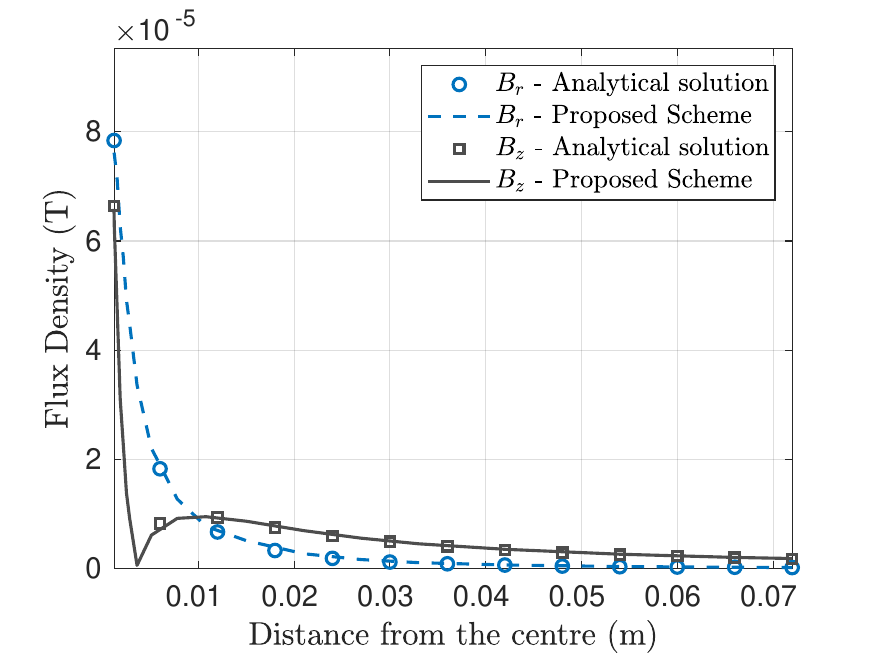}}}\\
		\mbox{\subfloat[]{\label{t9mur50} 
		\includegraphics[scale=0.62]{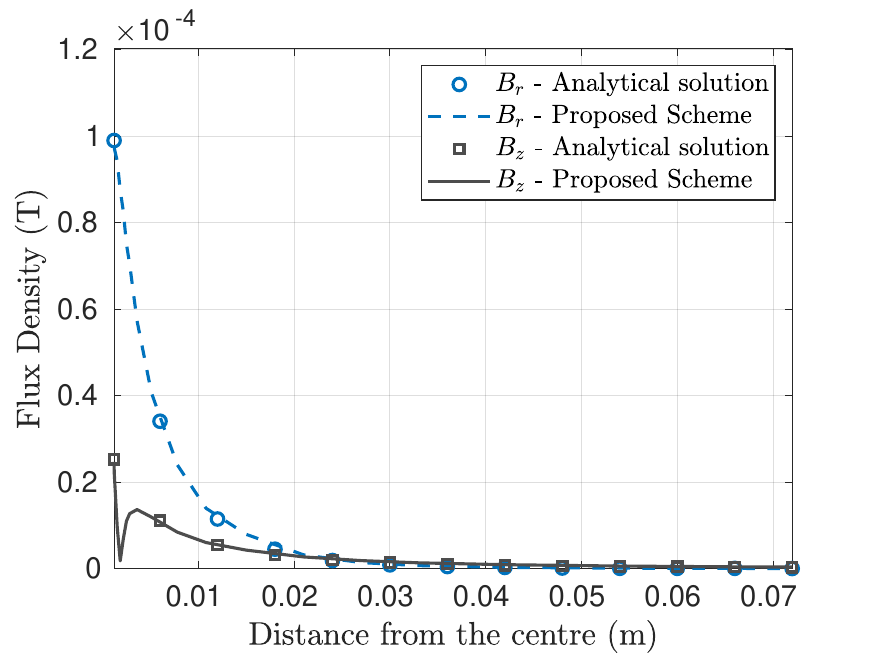}}}

		\caption{Comparison of the total magnetic flux densities from the
        analytical solution of the TEAM-9 problem \cite{team9} and the proposed
        formulation; for the cases of
        (a) $u_z=100ms^{-1}$, $\mu_r =  1$
        (b) $u_z=100ms^{-1}$, $\mu_r = 50$}
		\label{3dteam9results}
\end{figure}

The TEAM-9 test problem is also provided with the set of analytical solution
for comparison \cite{team9}. It may be noted that, the analytical solutions
are provided along the radius of $r=13mm$, which is $1mm$ away from both
the current carrying coil and the ferromagnetic cylinder. Since the
measurement point is very close to the cylinder, as well as, the circular
coil, it is necessary to model them as accurately as possible. 
Such a modeling is not feasible with linear edge elements in cartesian coordinate system. Therefore, the problem is transferred to the cylindrical
coordinate system for the accuracy study. In the cylindrical coordinate system, the results of Fig. \ref{Edge_A} and Fig. \ref{Elem_A} are once again observed.
In addition to this, edge elements can accurately represent the simulation
domain in cylindrical coordinate system. Simulations are carried out for
$v=100ms^{-1}$ with $\mu_r=1$, as well as, the ferromagnetic case of $\mu_r=50$.
The results from $r=13mm$ are plotted along with the analytical solution
in Fig. \ref{3dteam9results}. It can be seen that the proposed formulation performs consistently in 3D as well.

%-------------------------------------------------------------------------------
%     Part 6. Discussion on meshing
%-------------------------------------------------------------------------------
\section{Discussion on Mesh} \label{sec_mesh}

\begin{figure}
		\centering
		\mbox{\subfloat[]{\label{t9wedge_mesh} 
		\includegraphics[scale=0.23]{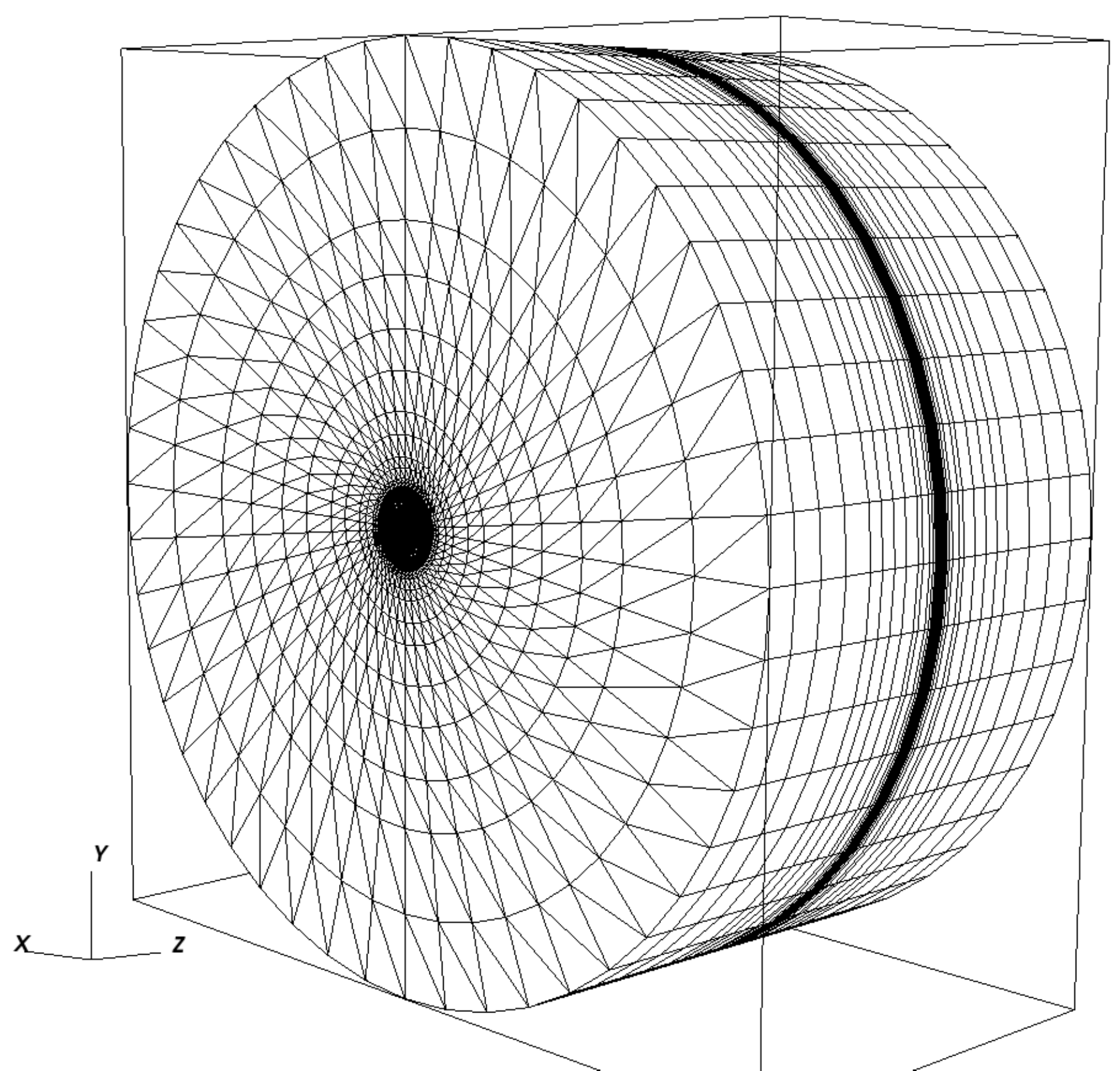}}}\\
		\mbox{\subfloat[]{\label{t9wedge_g} 
		\includegraphics[scale=0.55]{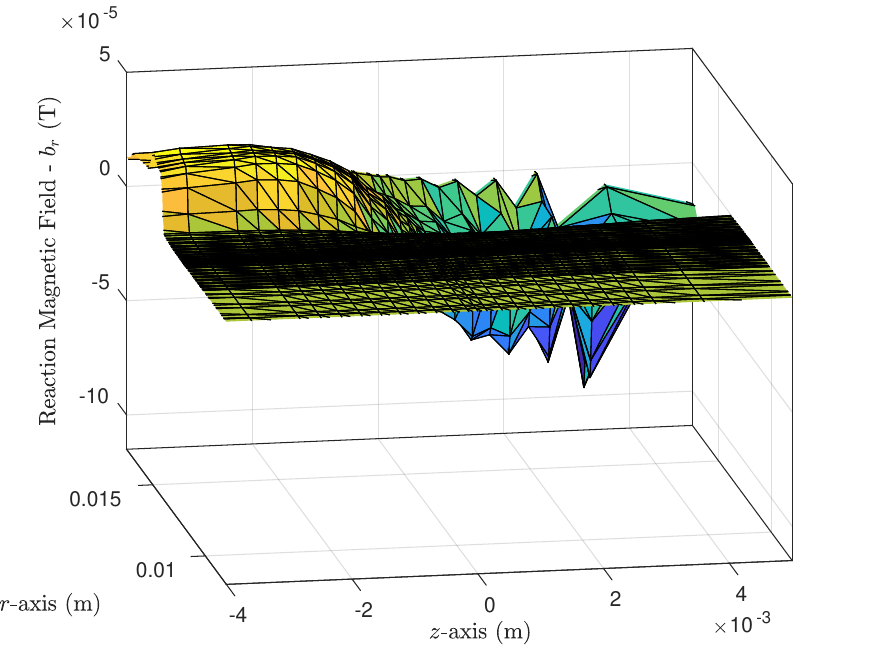}}}\\
		\mbox{\subfloat[]{\label{t9wedge_s} 
		\includegraphics[scale=0.55]{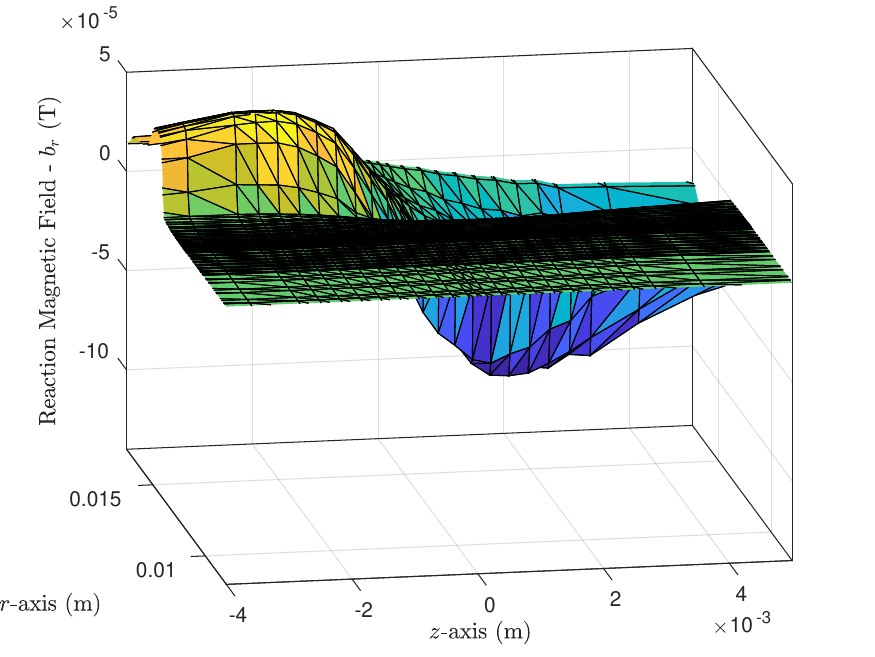}}}
		
		\caption{Sample results from TEAM 9a problem with wedge elements
		(a) Finite element mesh with wedge elements
		(b) Galerkin scheme - reaction magnetic field - $b_r$ for $u_z$ = 100 $ms^{-1}$ and $\mu_r$ = 50 
		(c) Proposed scheme - reaction magnetic field - $b_r$ for $u_z$ = 100 $ms^{-1}$ and $\mu_r$ = 50  in the cross section along $rz$ plane for 
$\theta \approx 0^o$.}
		\label{t9wedge}
\end{figure}

{The present work deals with the simulation of linear moving conductor 
problems such as, electromagnetic brakes, linear induction motor, 
electromagnetic flowmeter etc. 
In such cases, the conducting region of the problem can be 
and usually be discretised with graded regular mesh along the
moving direction.
In other words, the resulting mesh would look like a stack
of layers of different thickness along the moving direction.
The same can be seen in Fig. {\ref{team9a_3dmesh}}.
In this, the discretisation along the direction of motion ($z$-axis)
has dense discretisation close to the center, where
 the current loop is present and the discretisation becomes
coarser as we move away from the center.}

{The source based stabilisation strategies utilize this feature of the 
linear moving conductor problems.
Here, the stabilisation is brought in by the pole-zero cancellation of the
source term. Such an analysis is valid, only when the discretisation
is like a stack of layers along the moving direction. 
Hence, the proposed scheme requires a regular mesh in 2D with 
quadrilateral elements. In the case of 3D, the restriction only applies to the
direction of motion. Therefore, the cross section of the moving
conductor, i.e., the plane perpendicular to the motion, can be
discretised without any restrictions. 
Hence, the 3D problems can be discretised with hexahedral or wedge
elements. 
The discretisation using hexahedral elements is shown in 
Fig. {\ref{team9a_3dmesh}}.

The application of wedge elements with vector shape functions are scarce in
literature. However, the vector shape functions for the wedge elements are 
straight forward to derive and they are provided in Appendix {\ref{appA}} for 
quick reference. 
The discretisation using wedge elements for the TEAM-9 test problem 
is shown in Fig. {\ref{t9wedge_mesh}}. The simulated reaction magnetic field, 
along the $rz$-plane is plotted in Fig. {\ref{t9wedge_g}} and 
Fig. {\ref{t9wedge_s}}
for the Galerkin scheme and the proposed formulation, respectively. 
The $b_r$ from the proposed formulation is stable
as expected, with the wedge elements as well. 
As an added note, the discussion in this section is also applicable to 
node elements} \cite{su2}.
%The meshing restriction discussed here is common to both node and edge elements}

%-------------------------------------------------------------------------------
%     Part 7. Summary and Conclusion
%-------------------------------------------------------------------------------
\section{Summary and Conclusion}

Edge elements are vital in the finite element simulation of electromagnetic
fields, especially when multiple materials are present and the simulation
variables are electric and magnetic fields themselves. Like many other central-weighted numerical schemes, edge elements also produce numerically
oscillating solutions for the simulation of moving conductor problems at high
velocities. In such a situation, the usual strategy is to employ the upwinding
strategies, which in a way introduce extra diffusion to stabilize the solution
\cite{revfic1, edgeup2, r:edgeup3}. However, the upwinding schemes are known to be susceptible to transverse-boundary error at the material interfaces
\cite{soldreview1, nemosu, sus3, ss1}.

In this work, the source-based stabilization strategies which are 
proposed for the nodal formulation, are extended for the edge elements.
The formulation require a graded regular mesh along the direction of motion.
The stability of the proposed formulation is analytically studied in 1D, as 
well as, 2D with edge elements. Then numerical exercises are carried out
for the verification of the proposed formulation. The simulation results in 2D demonstrate that the formulation produces stable, accurate and converging
solutions. The 3D simulation is carried out with the TEAM-9 problem and 
stable solutions are observed. Comparing the analytical solutions of the TEAM-9 problem and the simulation results, accuracy of the proposed formulation
is demonstrated in 3D.

%-------------------------------------------------------------------------------
%     Part A. Appendices
%%-------------------------------------------------------------------------------
\appendices
\section{Wedge element - vector shape functions}\label{appA}
\begin{figure}
		\centering
		\includegraphics[scale=0.42]{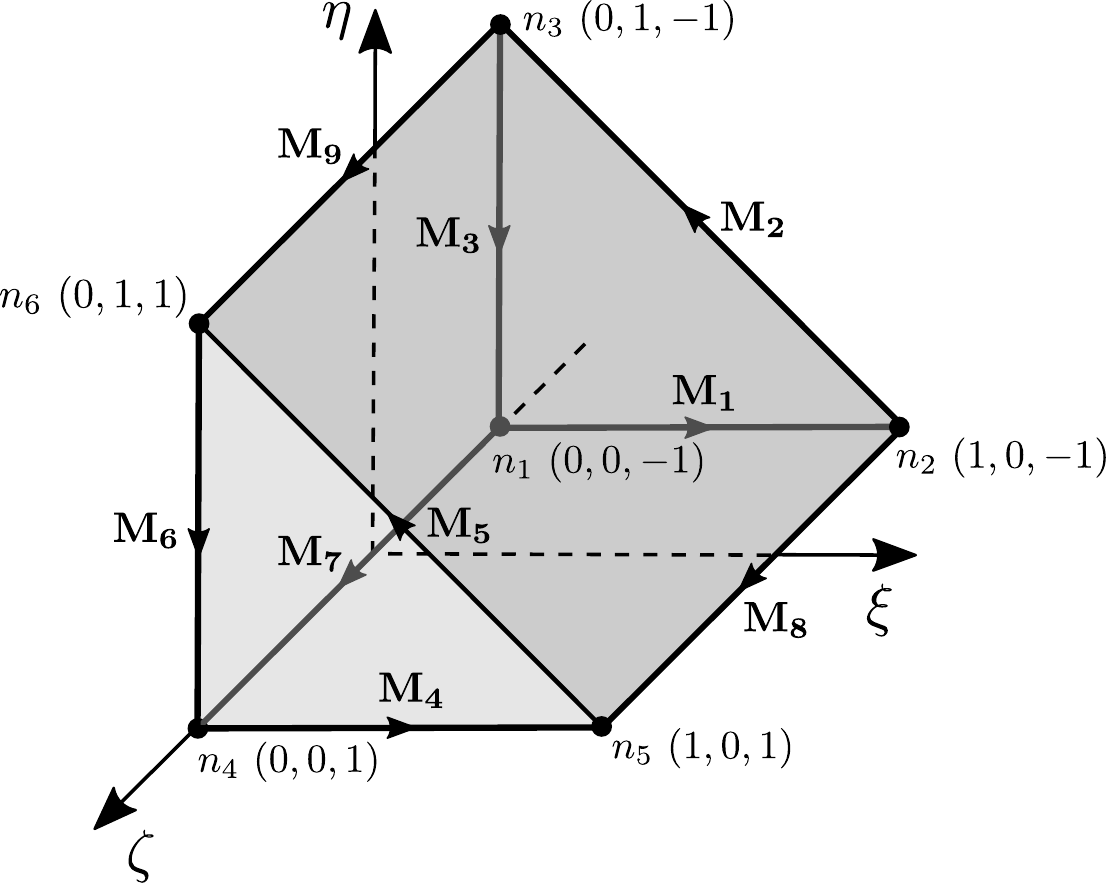}
		\caption{Wedge element in $(\xi, \eta, \zeta)$ coordinate system. 
		$n$-denotes the nodes of the element. $\mathbf{M}$-denotes the vector 
		shape functions, corresponds to each edge.}
		\label{wedge_elem}
\end{figure}

The edge shape functions for the wedge element are provided below. 
Consider the reference triangle element in $(\xi, \eta)$ coordinate system
with its nodes $a, b, c$ located at $(0,0)$, $(1,0)$, $(0,1)$ respectively.
The node shape functions for this reference triangle element can be written as \cite{febookr},
\begin{align}\label{eq:trishpfn}
N_a &= 1 - \xi -\eta \nonumber \\
N_b &=     \xi       \nonumber \\
N_c &=          \eta 
\end{align}
Using these, the vector shape function of the wedge element can be constructed.
The reference wedge element in $(\xi, \eta, \zeta)$ coordinate system is shown in Fig. {\ref{wedge_elem}}.
The edge shape functions for the wedge element are:
\begingroup
\allowdisplaybreaks
\begin{align}\label{eq:wedgeshpfn}
{\bf M_1} &= l_1(N_a ~\nabla N_b - N_b ~\nabla N_a)\dfrac{(1-\zeta)}{2}\nonumber \\
{\bf M_2} &= l_2(N_b ~\nabla N_c - N_c ~\nabla N_b)\dfrac{(1-\zeta)}{2}\nonumber \\
{\bf M_3} &= l_3(N_c ~\nabla N_a - N_a ~\nabla N_c)\dfrac{(1-\zeta)}{2}\nonumber \\
{\bf M_4} &= l_4(N_a ~\nabla N_b - N_b ~\nabla N_a)\dfrac{(1+\zeta)}{2}\nonumber \\
{\bf M_5} &= l_5(N_b ~\nabla N_c - N_c ~\nabla N_b)\dfrac{(1+\zeta)}{2}\nonumber \\
{\bf M_6} &= l_6(N_c ~\nabla N_a - N_a ~\nabla N_c)\dfrac{(1+\zeta)}{2}\nonumber \\
{\bf M_7} &= l_7 N_a ~\nabla \zeta /2 \nonumber \\
{\bf M_8} &= l_8 N_b ~\nabla \zeta /2 \nonumber \\
{\bf M_9} &= l_9 N_c ~\nabla \zeta /2
\end{align}
\endgroup
where, $l_1,~ l_2,~ \dots , ~l_9$ are the actual lengths of the edges that correspond
to edge shape function $\mathbf{M}$. The gradients are taken with respect to
the ($x, y, z$) coordinate system \cite{emfembookjin}.

%-------------------------------------------------------------------------------
%     Part R. References
%-------------------------------------------------------------------------------
% trigger a \newpage just before the given reference
% number - used to balance the columns on the last page
% adjust value as needed - may need to be readjusted if
% the document is modified later
%\IEEEtriggeratref{8}
% The "triggered" command can be changed if desired:
%\IEEEtriggercmd{\enlargethispage{-5in}}

% references section

% can use a bibliography generated by BibTeX as a .bbl file
% BibTeX documentation can be easily obtained at:
% http://mirror.ctan.org/biblio/bibtex/contrib/doc/
% The IEEEtran BibTeX style support page is at:
% http://www.michaelshell.org/tex/ieeetran/bibtex/
%\bibliographystyle{IEEEtran}
% argument is your BibTeX string definitions and bibliography database(s)
%\bibliography{IEEEabrv,../bib/paper}
%
% <OR> manually copy in the resultant .bbl file
% set second argument of \begin to the number of references
% (used to reserve space for the reference number labels box)
\bibliographystyle{IEEEtran} 
\bibliography{./References/References}
%\begin{thebibliography}{1}
%\bibitem{IEEEhowto:kopka}
%H.~Kopka and P.~W. Daly, \emph{A Guide to \LaTeX}, 3rd~ed.\hskip 1em plus
%  0.5em minus 0.4em\relax Harlow, England: Addison-Wesley, 1999.
%
%\end{thebibliography}

% biography section
% 
% If you have an EPS/PDF photo (graphicx package needed) extra braces are
% needed around the contents of the optional argument to biography to prevent
% the LaTeX parser from getting confused when it sees the complicated
% \includegraphics command within an optional argument. (You could create
% your own custom macro containing the \includegraphics command to make things
% simpler here.)

\begin{IEEEbiography}[{\includegraphics[width=1in,height=1.25in,clip,keepaspectratio]{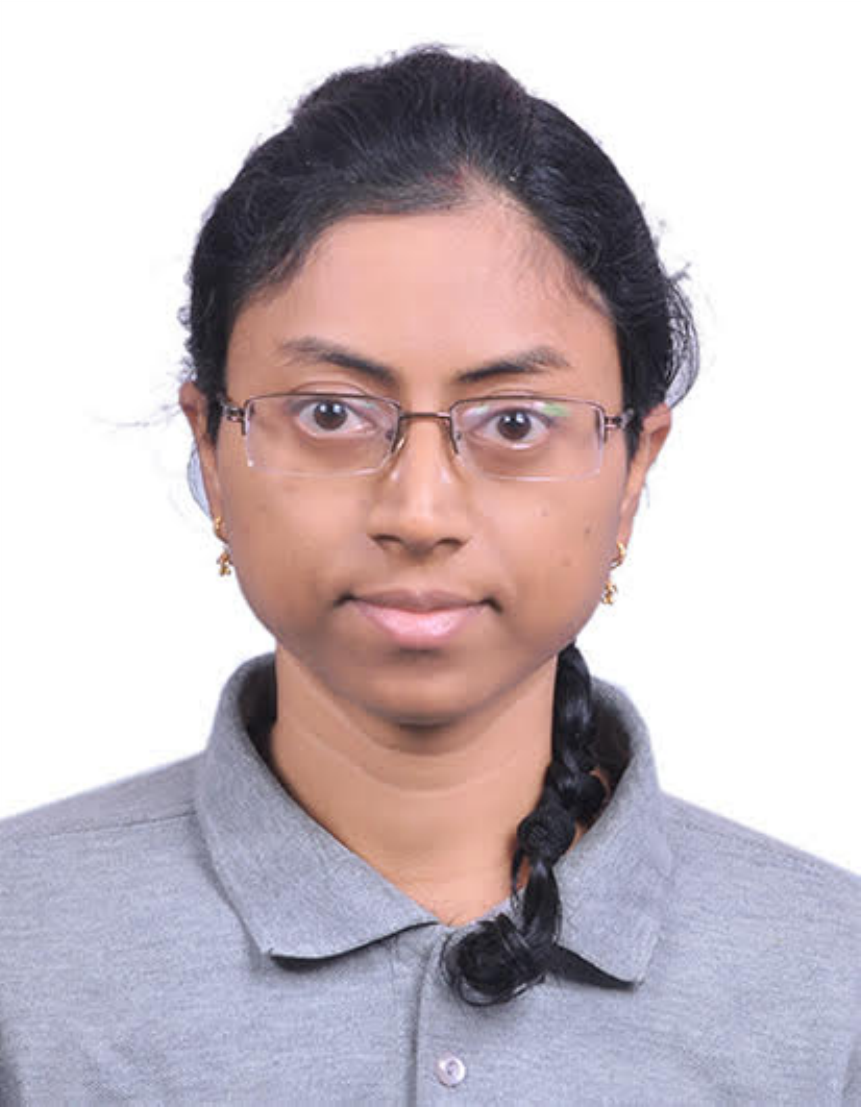}}]{Sujata Bhowmick}

 (Member, IEEE) received the B.E. degree in electrical engineering from IIEST, Shibpur, India, in 2006, and the M.E. degree in electrical engineering from the Indian Institute of Science, Bengaluru, India, in 2011. She received
the Ph.D. degree from the Department of Electronic Systems Engineering,
Indian Institute of Science, Bengaluru, India, in 2019.

Her current research interests include power electronics for renewable resources, single-phase grid-connected power converters, 
computational electromagnetics, finite element and 
edge element methods.

%Biography text here.
\end{IEEEbiography}

\begin{IEEEbiography}[{\includegraphics[width=1in,height=1.25in,clip,keepaspectratio]{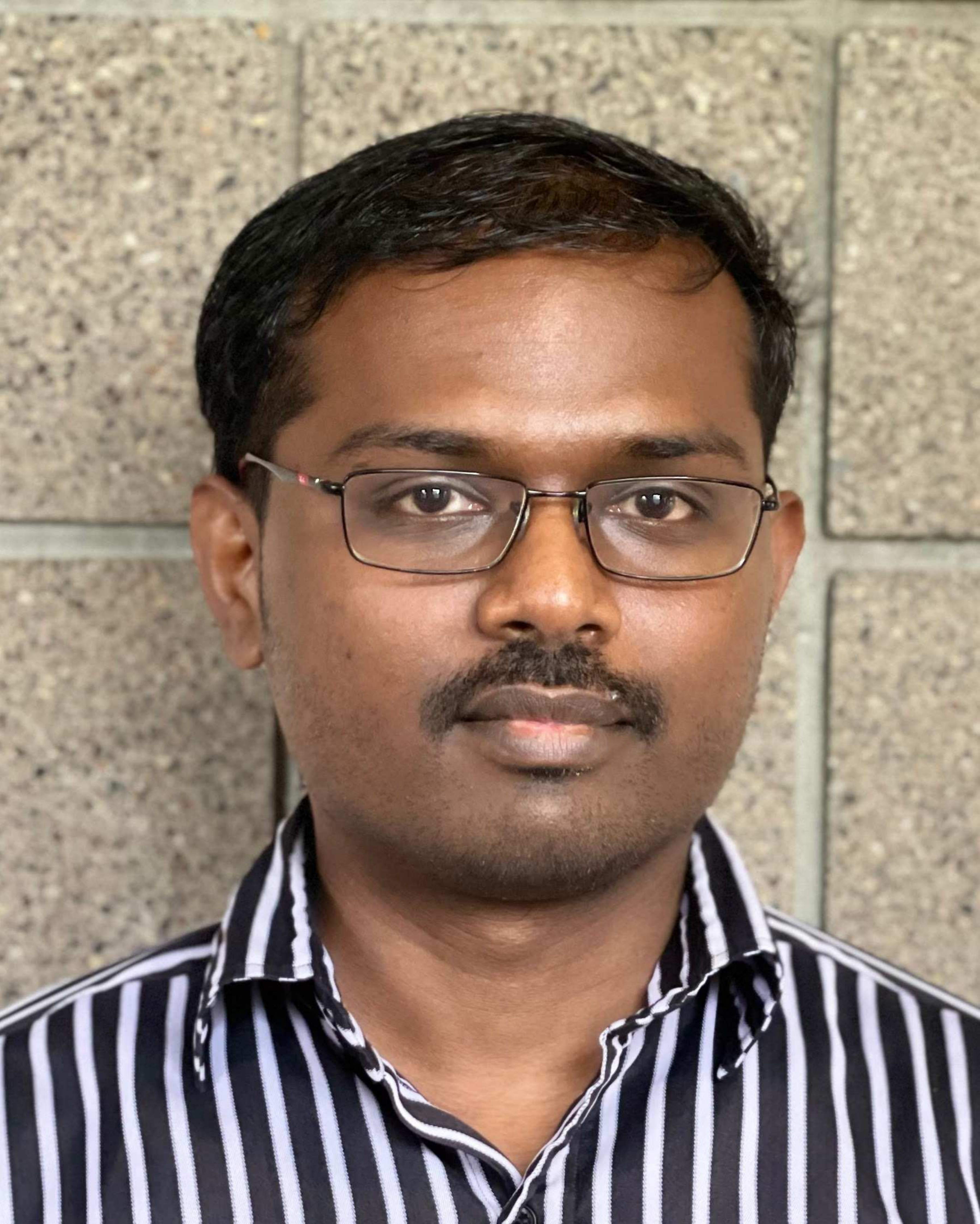}}]{Sethupathy Subramanian}
received the bachelors degree in electrical and electronics engineering from Anna University, Chennai, India, in 2009. He received the masters and doctrate degrees in electrical engineering from Indian Institute of Science, Bangalore, India in 2011 and 2017 respectively.
He also received masters and doctrate degrees in Physics from
University of Notre Dame, USA in 2022 and 2023 respectively.

His research interests, pertinent to electrical engineering, 
include computational electromagnetics, numerical stability, 
finite element and edge element methods. 

%Biography text here.
\end{IEEEbiography}
%

%\begin{IEEEbiography}[{\includegraphics[width=1in,height=1.25in,clip,keepaspectratio]{mshell}}]{Michael Shell}
% or if you just want to reserve a space for a photo:
%
%\begin{IEEEbiography}{Michael Shell}
%Biography text here.
%\end{IEEEbiography}
%
%% if you will not have a photo at all:
%\begin{IEEEbiographynophoto}{John Doe}
%Biography text here.
%\end{IEEEbiographynophoto}
%
%% insert where needed to balance the two columns on the last page with
%% biographies
%%\newpage
%
%\begin{IEEEbiographynophoto}{Jane Doe}
%Biography text here.
%\end{IEEEbiographynophoto}

% You can push biographies down or up by placing
% a \vfill before or after them. The appropriate
% use of \vfill depends on what kind of text is
% on the last page and whether or not the columns
% are being equalized.

%\vfill

% Can be used to pull up biographies so that the bottom of the last one
% is flush with the other column.
%\enlargethispage{-5in}

% that's all folks
\end{document}